\documentclass[journal ]{new-aiaa}
\usepackage[utf8]{inputenc}
\usepackage{textcomp}

\usepackage{graphicx}
\usepackage{amsmath}
\usepackage[version=4]{mhchem}
\usepackage{siunitx}
\usepackage{longtable,tabularx}
\usepackage{graphicx}
\usepackage{subfigure}
\usepackage[many]{tcolorbox}
\usepackage{algorithm} 
\usepackage{algpseudocode} 
\usepackage{mathtools}

\newtheorem{problem}{Problem}
\newtheorem{lemma}{Lemma}

\newtheorem{remark}{Remark}

\newcounter{example}[section]

\makeatletter
\usetikzlibrary{calc}
   
\makeatother

\title{Closed-Form Solutions for Minimum-Time Paths of Dubins Airplane in Steady Wind}

\author{Fanchen Wu \footnote{Ph.D student, Department of Aeronautics and Astronautics; fanchen-w@zju.edu.cn.}}
\affil{Zhejiang University, Hangzhou 310027, Zhejiang, China}

\author{Zheng Chen \footnote{Researcher, Department of Aeronautics and Astronautics;  z-chen@zju.edu.cn.}
}
\affil{Zhejiang University, Hangzhou 310027, Zhejiang, China}
\affil{Huanjiang Laboratory, Zhuji 311899, Zhejiang, China}

\begin{document}

\maketitle

\begin{abstract}
This paper is concerned with the minimum-time path-planning problem for a Dubins airplane under the influence of steady wind. The path-planning problem, by transforming into the air-relative frame, is equivalent to finding the minimum-time control strategy for a Dubins airplane to intercept a moving target. In the air-relative frame, by applying Pontryagin's maximum principle, the candidates for the minimum-time solution are categorized into a family of four types: $SC$, $CC$, $CCC$, $CSC$, where $S$ denotes a straight line segment and $C$ denotes a circular segment. Furthermore, the geometric properties for each type are analyzed, indicating that the paths of $SC$ and $CC$ can be obtained by finding the roots of a quadratic equation, while the paths of $CCC$ and $CSC$ are determined by the roots of some nonlinear transcendental equations. An improved bisection method is presented in the paper so that all the roots of the transcendental equations can be computed within a constant time. As a consequence, the globally optimal path can be obtained within a constant time by comparing all the candidates of the four types. Finally, numerical examples are presented, showing that the closed-form solutions derived in the paper can ensure to find the globally optimal solution by comparing with existing methods in the literature. 

\end{abstract}

% \iffalse
% \section*{Nomenclature}

% %\noindent(Nomenclature entries should have the units identified)

%  {\renewcommand\arraystretch{1.0}
%  \noindent\begin{longtable*}{@{}l @{\quad=\quad} l@{}}
 
%  $C$ & arc of unit circle \\
%  $d$ & length of straight line segment \\
%  $L$ & counterclockwise turn \\
%  $n$ & integer variable \\
%  $R$ & clockwise turn \\
%  $S$ & straight line segment \\
%  $s,t$ & binary set \\
%  $t_{f},T_{f}$ & terminal time \\
%  $u$ & control input \\
%  $V_T$ & target speed \\
%  $\boldsymbol{z}_0$ & initial state of the Dubins airplane \\
%  $[\bar{x}_{0},\bar{y}_{0}]$ & initial state of target \\
%  $[x_{f},y_{f}]$ & terminal position of the Dubins airplane\\
%  $[x_{T},y_{T}]$ & target position \\
%  $\alpha, \beta, \gamma$ & radian of arc \\
%  $\boldsymbol{c}_0^{{r}}, \boldsymbol{c}_0^{{l}}, \boldsymbol{c}_f^{{r}}, \boldsymbol{c}_f^{{l}}$ & center of unit circle \\
%  $\phi_{f}$ & terminal impact angle \\
%  $\tau_1, \tau_2$ & first and second switch time \\
%  $\theta_{P_0}$ & initial heading angle of the Dubins airplane \\
%  $\theta_{P_f}$ & terminal heading angle of the Dubins airplane \\
%  $\omega$ & wind speed \\

%  %\multicolumn{2}{@{}l}{Subscripts}\\
 
%  \end{longtable*}}

% \fi

\section{Introduction}

Fixed-wing Unmanned Aerial Vehicles (UAVs) are becoming increasingly popular in replacing manned aerial vehicles in many practical applications. Examples include tasks containing the element of danger, long endurance missions, and monotonous yet assiduous operations \cite{SHANMUGAVEL20101084}. The operation of fixed-wing UAVs without human intervention relies on advances in many domains, including automatic reasoning, perception, real-time control, {\it etc} \cite{https://doi.org/10.1049/tje2.12333}. One of the key topics in reasoning is path planning, which provides vehicles with the capability of deciding what motion commands to execute in order to achieve specified mission objectives \cite{8101650,1527001,yoon2019path}. Therefore, the foremost feature of path planners for fixed-wing UAVs, especially in highly dynamic scenarios with time-critical missions, is the ability to generate optimal paths {\it in situ} while satisfying various constraints. 

In realizing real-time path planning for fixed-wing UAVs, there is an important trade-off between utilizing a realistic model and tractability. Balancing this trade-off, the model of unidirectional non-holonomic vehicles, moving at a constant speed with a minimum turning radius, was conceived for path planning of fixed-wing UAVs \cite{mclain2014implementing,6842272,2018Time}. Following the work \cite{markov1887some} by A. A. Markov in 1887 and the work \cite{10.2307/2372560} by L. E. Dubins in 1957, such non-holonomic vehicles are usually dubbed as Markov-Dubins vehicles or simply Dubins vehicles. This vehicle model provides a highly accurate approximation to the kinematics of fixed-wing UAVs in altitude-hold mode. Therefore, some well-established results for the minimum-time paths of Dubins vehicles have been directly applied to planning paths for fixed-wing UAVs \cite{MATVEEV2011515,doi:10.2514/1.44580,doi:10.2514/1.G005748,CHEN2020108996}. In the literature, a fixed-wing UAV, moving at an altitude-hold mode with a constant speed, is also called Dubins airplane; see, e.g., \cite{2008Time,choi2014time}.

It should be noted that the minimum-time path is equivalent to the shortest path as the speed is considered constant. According to \cite{10.2307/2372560}, the shortest path for a Dubins airplane between two configurations (including a 2-dimensional position and a heading angle) is a $C^1$ path, which is a smooth concatenation of circular arcs and straight line segments. Further studies by geometric analysis in \cite{10.2307/2372560} and by optimal control in \cite{Sussmann1991SHORTESTPF,351019} showed that the shortest path lies in a sufficient family of 6 types. This allows computing the shortest path within a constant time by checking at most 6 solution candidates. 
Up to present, the shortest paths for Dubins airplanes in more complex scenarios have been studied. For example, the shortest path from a configuration to a target circle was studied in \cite{doi:10.2514/6.2019-0919,CHEN2020108996}, and analytical solutions were presented by using the Pontryagin's Maximum Principle (PMP). The shortest path for Dubins airplanes between two configurations via an intermediate point was studied in \cite{CHEN2019368}; by restricting the shortest path into a family of 18 types, a polynomial method was proposed to compute the shortest path. A more complex scenario involving obstacle avoidance was further considered in \cite{JONES2021109510,JHA2022110637}. In \cite{JONES2021109510}, a discretization-based dynamic programming framework was proposed for planning optimal paths while avoiding static obstacles. To incorporate the avoidance of moving obstacles, Jha {\it et al.} in \cite{JHA2022110637} obtained the minimum-time path by solving a set of nonlinear equations.

Although the papers cited in the preceding paragraph provide some approaches for path planning of Dubins airplanes, they cannot be applied to the scenario with the effect of wind. 
In recent decades, the minimum-time paths for Dubins airplanes in steady wind have been extensively studied.
Rysdyk {\it et al.} in \cite{doi:10.2514/1.27359} observed that the path for a Dubins airplane flying in steady wind with a constant maximum turn rate can be formulated as a trochoidal curve. Techy {\it et al.} in \cite{4739456} further found that the minimum-time path in steady wind consists only of straight-line segments and trochoidal segments. 
This allows for the generation of feasible paths by simply concatenating straight-line segments and trochoidal segments.
In a subsequent work \cite{doi:10.2514/1.44580}, the authors categorized the minimum-time path into six types. Among the six types, two types could be analytically computed as closed-form solution were obtained, but the other four should be computed using numerical root-finding algorithms to solve some transcendental equations. It is widely known that employing numerical methods often require significant computational efforts, restricting its applications in real-time scenarios. Additionally, if the initial guess is not appropriately chosen for transcendental equations, the solutions may fall into local optima or even fail to converge \cite{doi:10.2514/3.19932}.

The studies mentioned above were conducted in inertial frames, where the path was distorted by the influence of wind. To avoid these distortions, 
it is found that the minimum-time path of a Dubins airplane affected by steady wind in the air-relative frame, an inertial frame moving in the direction of the wind, consists of straight-line segments and circular arcs with maximum turn rate \cite{doi:10.2514/6.2005-6186}. In addition, the problem of minimum-time path planning in steady wind is equivalent to a moving-target interception problem, where the virtual target moves with the same speed as the steady wind but in the opposite direction \cite{bakolas2013optimal}. For this reason, the minimum-time path of a Dubins airplane in steady wind have been studied in the literature by analysing the solution path of moving-target interception problem; see, e.g., \cite{doi:10.2514/6.2005-6186, mcgee2007optimal,doi:10.2514/1.G001527, gopalan2017generalized,MITTAL2020103646,ZHENG2021}. 
In \cite{doi:10.2514/6.2005-6186, mcgee2007optimal}, the optimal solutions for the moving-target interception problem were categorized into a family of ten types, and numerical methods were presented for each type.
To be able to compute the solution path in real time,
Gopalan {\it et al.} \cite{doi:10.2514/1.G001527, gopalan2017generalized} explored the closed-form solutions for four of the ten types, under the assumption that the initial position of the Dubins airplane is at least four times the minimum turning radius away from the target; closed-form solutions were derived for two types, while the remaining two required finding the zeros of a nonlinear equation by numerical methods.
The two types of paths with closed-form solutions were further studied by Mittal {\it et al.} in \cite{MITTAL2020103646}, showing that full reachability of any target, provided it is at least four times the minimum turning radius from the initial position of the Dubins airplane, can be achieved by extending the ranges of the arc segments for these two types of paths.
Recently, by analysing the geometric properties of the four path types presented in \cite{doi:10.2514/1.G001527}, Zheng {\it et al.} in \cite{ZHENG2021} transformed the minimum-time path-planning problem into the problem of finding the roots for certain nonlinear equations and proposed a robust method to determine all the roots of these equations.
In \cite{doi:10.2514/1.G001527,gopalan2017generalized,MITTAL2020103646,ZHENG2021}, all closed-form solutions were derived under the assumption of a sufficient initial distance between the Dubins airplane and the target. However, to the authors' best knowledge, research on closed-form solutions for minimum-time paths in steady wind without any assumption on the distance between initial position and target position is scarce to see. 

In this paper, we present a more systematic and efficient method for obtaining closed-form solutions to the minimum-time path-planning problem in steady wind.
To achieve this, the path-planning problem is formulated as an optimal control problem in air-relative frame. 
Using the PMP, the solution paths for such an optimal control problem are categorized into a sufficient family of four types. To derive closed-form solutions, the geometric properties are analysed for each type. While explicit expressions are obtainable for the paths of two types, the remaining two require solving a series of transcendental equations. Given that each transcendental equation may have multiple roots and existing numerical methods could converge to local optima, an improved bisection method is proposed to find all the roots of these equations.
To this end, one is allowed to obtain the globally minimum-time solution by comparing all the candidate types.

The paper is organized as follows: 
Section \ref{SE:preliminary} formulates the minimum-time path planning problem in steady wind and presents the equivalent optimal control problem in the air-relative frame. Necessary conditions for optimality are presented by applying the PMP and the solution paths are categorized into a sufficient family of four types in Section \ref{SE:Characterization}. By analysing the geometric properties of each type, closed-form solutions are provided for each type in Section \ref{SE:Analytical}.
Numerical simulations are presented in Section \ref{SE:simulation}, verifying the developments of this paper.

\section{Problem Formulation}\label{SE:preliminary}

\begin{figure}[hbt!]
    \centering
    \includegraphics[width=.5\textwidth]{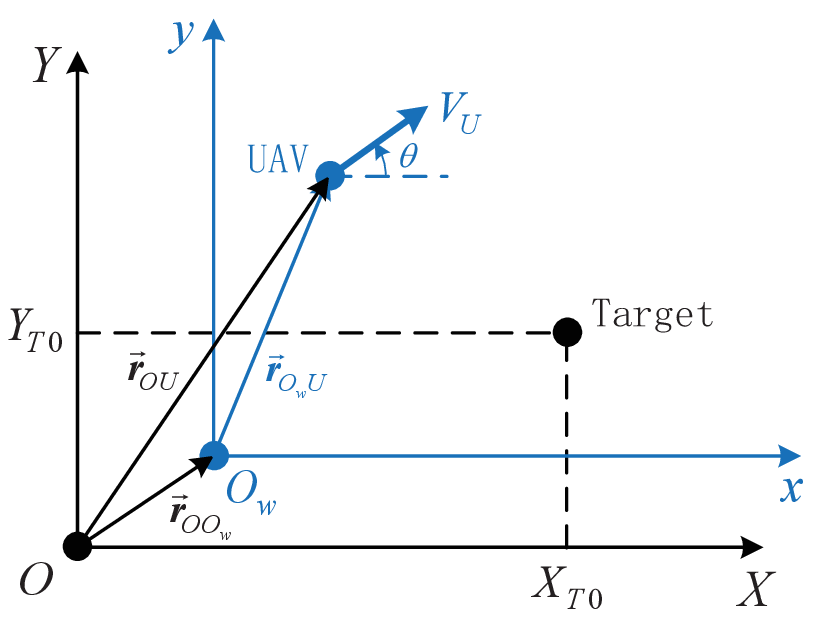}
    \caption{Geometry of two reference frames.}
    \label{Fig:geometry1}
\end{figure}

Consider a scenario that a Dubins airplane moves in altitude-hold mode with a constant airspeed, and we aim to solve a minimum-time control problem of steering the Dubins airplane, affected by a steady wind, to a desired target with an expected final heading angle. In order to formulate the optimal control problem, we use an air-relative frame to describe the motion of the Dubins airplane. We first consider an inertial Cartesian coordinate frame $OXY$ with its origin being the same as the initial position of the Dubins airplane. The $OX$ axis points to the East, and the $OY$ axis points to the North. Correspondingly, the air-relative frame $O_w xy$, which moves with the steady wind, is introduced as illustrated in Fig.~\ref{Fig:geometry1}. The origin of frame $O_w xy$ relative to frame $OXY$ is denoted by the vector ${\mathop{\pmb{r}}\limits ^{\rightarrow}}_{OO_w}$, the x-axis points to the East, and the y-axis points to the North.

Denote by $(x,y)\in \mathbb{R}^2$ the coordinate position of the Dubins airplane in the air-relative frame, and denote by $\theta\in [0,2\pi)$ the corresponding heading angle, which is measured counterclockwise from the $x$-axis to the velocity vector of the airplane. Then, in the air-relative frame, if we normalize the speed to one \cite{doi:10.2514/6.2005-6186}, the kinematics for the Dubins airplane in altitude-hold mode can be expressed as 
\begin{equation}
    \label{Eq:system2}
    \left\{
        \begin{aligned}
		\dot{x}(t) = & \cos \theta(t) \\
		\dot{y}(t) = & \sin \theta(t) \\
		\dot{\theta}(t) = & \dfrac{u(t)}{\rho}
		\end{aligned}
	\right.
\end{equation}
where $t\geq 0$ denotes the time, the over dot denotes the differentiation with respect to time, $\rho$ is a positive constant related to the minimum turning radius of the Dubins airplane, and $u\in [-1,1]$ is the control parameter, representing the lateral acceleration. 

Denote by $(w_x,w_y)\in \mathbb{R}^2$ the velocity of the wind in the inertial frame $OXY$. Then, the minimum-time control problem of steering the Dubins airplane to a desired target with an expected heading angle, if stated in the air-relative frame $O_wxy$, is equivalent to the following optimal control problem.
\begin{problem}
\label{problem2}
Find a minimum time $t_f > 0$ and a measurable control $u(\cdot)\in[-1,1] $ over the interval $[0,t_f]$ so that the system in Eq.~(\ref{Eq:system2}) is steered from a fixed initial state $(x_0,y_0,\theta_0)$ to a target $(x_f,y_f,\theta_f)$ satisfying 
$$(x_f,y_f) = (X_{T0},Y_{T0}) - t_f (w_x,w_y)$$
where $(X_{T0},Y_{T0})$ is the position of the target in the inertial frame $OXY$.
\end{problem}
Throughout the paper, we assume that the speed of the wind is smaller than one, i.e., $w_x^2 + w_y^2 < 1$, to make sure the existence of solutions for Problem \ref{problem2} \cite{doi:10.1137/0115133}. Without loss of generality, we fix the initial state of the Dubins airplane at $t=0$ in the air-relative frame as $(0,0,\pi/2)$, i.e.,
$$(x_0,y_0,\theta_0) =(0,0,\pi/2)$$

The optimal control problem in Problem \ref{problem2} has been extensively studied in the literature; see, e.g., \cite{ mcgee2007optimal,doi:10.2514/1.G001527, gopalan2017generalized,MITTAL2020103646,ZHENG2021}, because it is fundamentally important for planning minimum-time paths for fixed-wing UAVs in steady wind and for Unmanned Surface Vehicles (USVs) with constant ocean currents  \cite{hao2023optimal}. In the next two sections, the solutions to Problem \ref{problem2} will be classified into four categories by using the PMP, and the geometric properties of each category will be analysed to derive closed-form solutions. 

\section{Characterization of the Optimal Solution}\label{SE:Characterization}

Denote by $\pmb{p}=(p_x,p_y,p_\theta) \in \mathbb{R}^3$ the costate vector of $(x,y,\theta)$. Then, the Hamiltonian for Problem \ref{problem2} is expressed as
\begin{align}
\label{Eq:H_aaa}
	H = p_x \cos\theta + p_y \sin\theta + p_\theta\dfrac{u}{\rho} - 1
	%\nonumber
\end{align}
The costate variables are governed by
\begin{align}
	\dot{p}_x(t) &= -\dfrac{\partial H}{\partial x}=0 \label{Eq:px}\\
	\dot{p}_y(t) &= -\dfrac{\partial H}{\partial y}=0 \label{Eq:py}\\
	\dot{p}_\theta(t) &= -\dfrac{\partial H}{\partial \theta}= p_x(t) \sin\theta(t) - p_y(t) \cos\theta(t)  \label{Eq:ptheta}
\end{align}
It is apparent from Eq.~(\ref{Eq:px}) and Eq.~(\ref{Eq:py}) that $p_x$ and $p_y$ are constant along an optimal path. By integrating Eq.~(\ref{Eq:ptheta}), we have
\begin{align}
	p_\theta(t) =  p_x y(t) - p_y x(t) + c_0
 \label{Eq:ptheta2}
\end{align}
where $c_0$ is a constant. In view of Eq.~(\ref{Eq:ptheta2}), if $p_\theta$ is identical to zero on a nonzero interval, the path is a straight line segment. Thus, we have 
\begin{align}
\label{Eq:u_0_ptheta}
    u\equiv0,\quad \rm{if} \ p_\theta\equiv 0
\end{align}
According to Eq.~(\ref{Eq:u_0_ptheta}) and the PMP \cite{1987Pontryagin}, the optimal control $u$ is totally determined by $p_\theta$, i.e.,
\begin{equation}
    \label{Eq:u}
    u=\left\{
        \begin{aligned}
		&1, \quad &p_\theta>0\\
		&0, \quad &p_\theta=0\\
		&-1, \quad &p_\theta<0
		\end{aligned}
	\right.
\end{equation}
The path is a circular arc turning clockwise (resp. counterclockwise) if $u=-1$ (resp. $u=1$), or a straight line segment if $u=0$. Therefore, the switching conditions in Eq.~(\ref{Eq:u}) imply that the solution for Problem \ref{problem2} is a concatenation of circular arcs (denoted by "$C$") and straight line segments (denoted by "$S$"). If the arc "$C$" turns clockwise (resp. counterclockwise), we specify it as "$R$" (resp. "$L$").
Additionally, for any $\alpha \in [0,2\pi)$, we denote by $C_{<\alpha}$, $C_{\alpha}$, and $C_{>\alpha}$ the corresponding circular arc with a radian smaller than, equal to, and greater than $\alpha$, respectively.

By the following remark, we shall show that the solution pattern of Problem \ref{problem2} belongs to four categories:
\begin{remark}[Buzikov and Galyaev \cite{BUZIKOV2022109968}]
The solution pattern of Problem \ref{problem2} belongs to the family $$\mathcal{F} \coloneqq \{SC_{2\pi},CC_{2\pi},CCC,CSC\}$$
and subpatterns in $\mathcal{F}$.
\label{Remark:1}
\end{remark}
Denote by $T_{SC}$ the time duration for the Dubins airplane to travel along the solution of type $SC_{2\pi}$; if the path of $SC_{2\pi}$ is not a feasible path for Problem \ref{problem2}, we set $T_{SC}$ as $+\infty$. The same definitions apply to $T_{CC}$, $T_{CCC}$, and $T_{CSC}$. Then, according to Remark \ref{Remark:1}, we have 
\begin{equation}
\label{Eq:t_f}
    t_f = \min \{T_{SC}, T_{CC}, T_{CCC}, T_{CSC}\}
\end{equation}
As a result of Eq.~(\ref{Eq:t_f}), it is apparent that the global minimum time for Problem \ref{problem2} can be obtained immediately once all the four types of paths can be found in real time. Some algebraic transcendental equations were presented in \cite{BUZIKOV2022109968} so that the solution path of each type can be computed by numerically finding the roots of a transcendental equation. However, a transcendental equation may have multiple roots, but the solution path is only related to a specific root, as shown by the numerical examples in Section \ref{SE:simulation}. Existing numerical algorithm cannot guarantee to find the desired root. Even if there is only one root for the transcendental equation, existing numerical algorithms may suffer the issue of convergence. Therefore, it is challenging to compute the paths of the four types via numerically finding roots of transcendental equations. To address this issue, in the next section so that the globally optimal solution of Problem \ref{problem2} can be computed within a constant time. 

\section{Closed-form Solutions}
\label{SE:Analytical} 

In this section, we will present the closed-form solution for each type in $\mathcal{F}$ by analysing its geometric property.

\subsection{Closed-Form Solution for Path of \texorpdfstring{$SC_{2\pi}$}{} }

\begin{figure}[htbp]
    \centering
    
    \subfigure[$SR_{2\pi}$]{
        \begin{minipage}[t]{7cm}
        \centering
        \includegraphics[height=4.3cm]{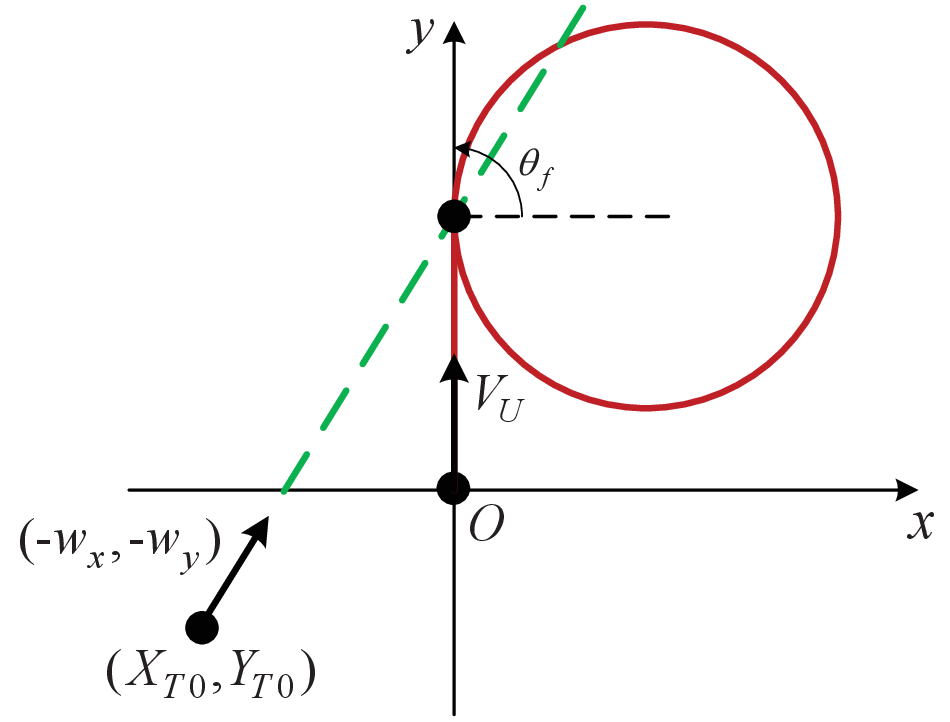}
        \label{Fig:geometry_SR2pi}
    \end{minipage}%
    }%
    \subfigure[$SL_{2\pi}$]{
        \begin{minipage}[t]{7cm}
        \centering
        \includegraphics[height=4.3cm]{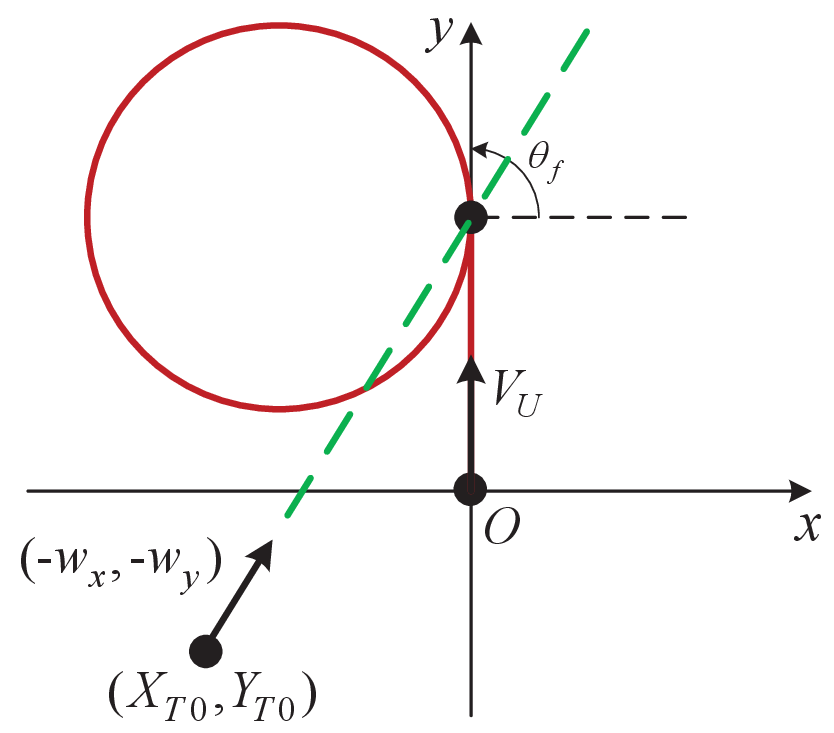}
        \label{Fig:geometry_SL2pi}
    \end{minipage}%
    }%
    
    \centering
    \caption{Geometry for the path of $SC_{2\pi}$ in the air-relative frame.}
    \label{Fig:geometrySC}
\end{figure}

The $SC_{2\pi}$-path includes two candidates, i.e.,
$$SC_{2\pi}=\{SR_{2\pi},SL_{2\pi}\}$$
both of which are illustrated in Fig.~\ref{Fig:geometrySC}. It is evident that to reach the same target, the $SR_{2\pi}$-path is symmetric to the $SL_{2\pi}$-path with respect to the $y$-axis, and the time durations for the Dubins airplane to travel along these two paths are identical. We first present the closed-form solution for the $SR_{2\pi}$-path.
%Therefore, in this subsection, we only present the closed-form solution for the $SR_{2\pi}$-path.

Denote by $U_{SR}$ the control along the $SR_{2\pi}$-path. Then, according to Eq.~(\ref{Eq:u}), we have 
\begin{equation}
    U_{SR}(t)=\left\{
    \begin{aligned}
        &0,\quad \quad &&t\in[0,d)\\
        &-1, \quad \quad &&t\in[d,d+2\pi\rho]
    \end{aligned}
    \right.
    \label{Eq:u_SR}
\end{equation}
where $d$ is the length of the straight-line segment of the $SR_{2\pi}$-path. Note that $d$ is equivalent to the time duration for the Dubins airplane to move along the straight-line segment since the speed is normalized to one. It is now enough to find the value of $d$ in order to find the solution path of $SR_{2\pi}$. 

According to Fig.~\ref{Fig:geometry_SR2pi}, we have that the solution path is of $SR_{2\pi}$ if and only if the following four equations are satisfied:
\begin{equation}
\left\{
    \begin{aligned}
    &x_f=0\\
    &y_f=d\\
    &\theta_f=\pi/2\\
    &-w_y(x_f-X_{T0})=-w_x(y_f-Y_{T0}) 
    \end{aligned}
\right.
\label{EQ:SC0}
\end{equation}
Note that the time required for the Dubins airplane to travel from its initial state $(x_0,y_0,\theta_0)$ to the desired terminal position $(x_f,y_f)$ is equal to the time required for the target to reach the same terminal position in the air-relative frame, indicating
\begin{align}
	d + 2\pi\rho= \dfrac{\sqrt{(x_f - X_{T0})^2 + (y_f - Y_{T0})^2}}{\sqrt{(-w_x)^2 + (-w_y)^2}}
	\label{EQ:SC2}
\end{align}
Combining Eq.~(\ref{EQ:SC2}) and Eq.~(\ref{EQ:SC0}), we have
\begin{align}
    -(d + 2\pi\rho)w_y = d - Y_{T0}
    \label{EQ:SC3}
\end{align}
Rearranging Eq.~(\ref{EQ:SC3}), the explicit expression for $d$ is given as
\begin{align}
    d=\dfrac{Y_{T0}-2\pi\rho w_y}{1+w_y}
    \label{EQ:SC4}
\end{align}
Substituting Eq.~(\ref{EQ:SC4}) into Eq.~(\ref{EQ:SC0}), it follows that the $SR_{2\pi}$-path exists if the following two equations hold:
\begin{equation}
\left\{
    \begin{aligned}
    &\theta_f=\pi/2\\
    &w_y X_{T0}-w_y Y_{T0}=\dfrac{2\pi\rho w_y-Y_{T0}}{1+w_y}
    \end{aligned}
\right.
\label{EQ:SC5}
\end{equation}
Then, by the definition of $T_{SC}$, we have
\begin{equation}\label{Eq:T_SR}
T_{SR} = \left\{
\begin{aligned}
    &d + 2\pi \rho, \quad && \text{if Eq.~(\ref{EQ:SC5}) is satisfied}\\
    &+\infty, \quad &&\text{if Eq.~(\ref{EQ:SC5}) is not satisfied}
\end{aligned}
\right.
\end{equation}
In addition, the control along the path of $SR_{2\pi}$ is determined by Eq.~(\ref{Eq:u_SR}) if Eq.~(\ref{EQ:SC5}) is satisfied. 

According to Fig.~\ref{Fig:geometrySC}, if the path type is $SL_{2\pi}$, $T_{SL}$ will take the same form as Eq.~(\ref{Eq:T_SR}), but the control along the path of $SL_{2\pi}$ is 
\begin{equation}
    U_{SL}(t)=\left\{
    \begin{aligned}
        &0,\quad \quad &&t\in[0,d)\\
        &1, \quad \quad &&t\in[d,d+2\pi\rho]
    \end{aligned}
    \right.
\end{equation}

\subsection{Closed-Form Solution for Path of \texorpdfstring{$CC_{2\pi}$}{} }

\begin{figure}[htbp]
    \centering
    
    \subfigure[$RL_{2\pi}$]{
        \begin{minipage}[t]{7cm}
        \centering
        \includegraphics[height=4.4cm]{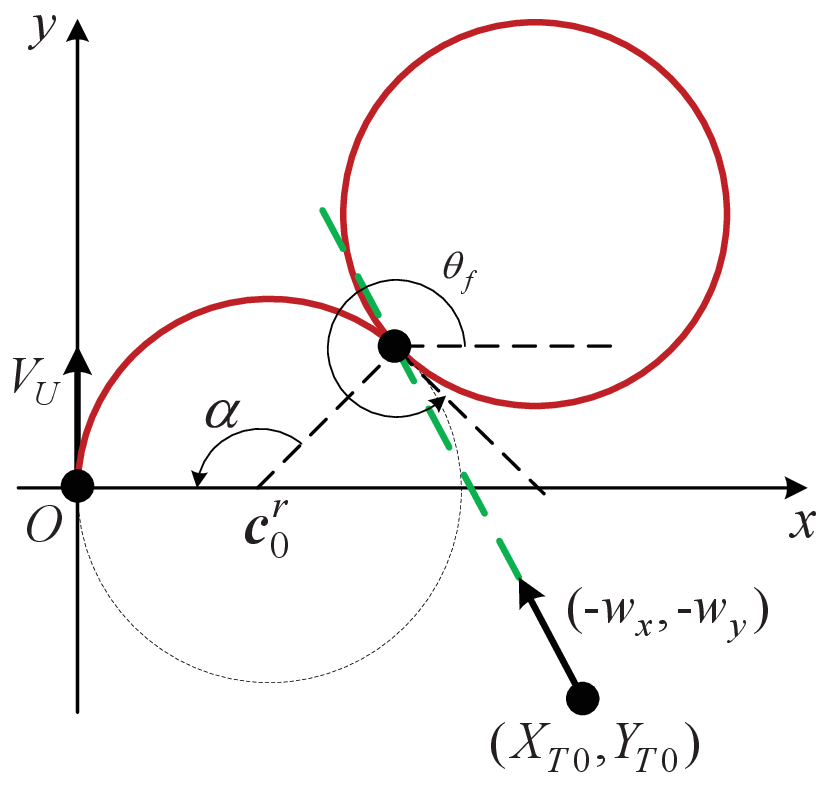}
        \label{Fig:geometry_RL2pi}
    \end{minipage}%
    }%
    \subfigure[$LR_{2\pi}$]{
        \begin{minipage}[t]{7cm}
        \centering
        \includegraphics[height=4.6cm]{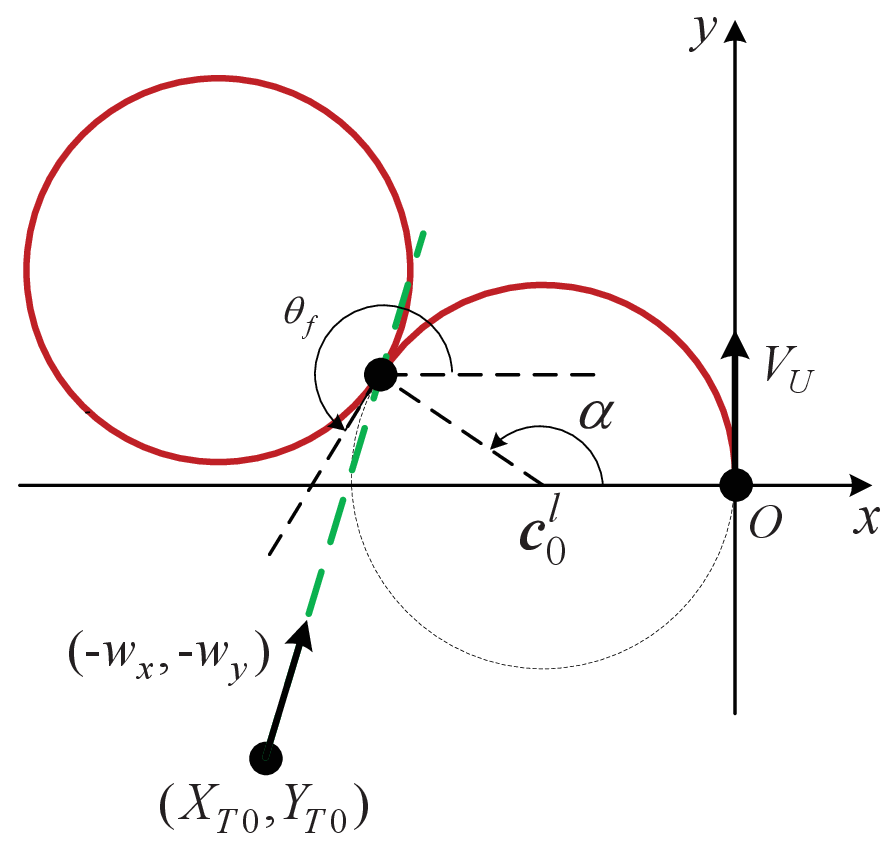}
        \label{Fig:geometry_LR2pi}
    \end{minipage}%
    }%
    
    \centering
    \caption{Geometry for the path of $CC_{2\pi}$ in the air-relative frame.}
    \label{Fig:geometryCC}
\end{figure}

The $CC_{2\pi}$-path includes two candidates, i.e.,
$$CC_{2\pi}=\{RL_{2\pi},LR_{2\pi}\}$$
Let $\alpha \in [0,2\pi)$ be the radian of the first circular arc, as shown in Fig.~\ref{Fig:geometryCC}. 
Denote by $U_{CC}$ the control along the $CC_{2\pi}$-path. Then, we have that $U_{CC}$ is expressed as
\begin{equation}\label{Eq:U_CC}
    U_{CC}(t)=\left\{
    \begin{aligned}
        &\sigma,\quad \quad &&t\in[0,\rho\alpha)\\
        &-\sigma, \quad \quad &&t\in[\rho\alpha,\rho(\alpha+2\pi)]
    \end{aligned}
    \right.
\end{equation}
where $\sigma=-1$ if the solution path is of $RL_{2\pi}$, and $\sigma=1$ if the solution path is of $LR_{2\pi}$. It is now enough to find the value of $\alpha$ in order to find the solution path of $CC_{2\pi}$. 
By the following lemma, we shall show that $\alpha$ is a root of a quadratic function.

\begin{lemma}
    \label{Le:CC}
    If the solution path for Problem \ref{problem2} is of type $CC_{2\pi}$, the following two statements hold:
    \begin{description}
    \item[(1)] The solution path is of $RL_{2\pi}$ if and only if 
    \begin{equation}\label{EQ:CC1_1}
    \left\{
        \begin{aligned}
        \theta_f &=\pi/2-\alpha+2n\pi\\
        w_x &= \dfrac {\dfrac{X_{T0}}{\rho}-\cos(\dfrac{\pi}{2}+\theta_f)-1}{\dfrac{\pi}{2}-\theta_f+2(n+1)\pi}\\
        w_y &= \dfrac {\dfrac{Y_{T0}}{\rho}-\sin(\dfrac{\pi}{2}+\theta_f)}{\dfrac{\pi}{2}-\theta_f+2(n+1)\pi}
        \end{aligned}
    \right.
    \end{equation}
    where $n=0$ (resp. $=1$) if $\theta_{f}<\pi/2$ (resp. $\geq \pi/2$); in addition, if Eq.~(\ref{EQ:CC1_1}) holds, we have that $\alpha$ is a root of 
    \begin{align} \label{EQ:CC1_2}
     a_1 \alpha^2+a_2\alpha+a_3=0   
    \end{align}
    where the constants $a_1$--$a_3$ are presented in Appendix \ref{Appendix:A}.
    % \begin{equation}
    % \left\{
    %     \begin{aligned}
    %     a_1&=\rho^2 w_y^2+\rho^2 w_x^2\\
    %     a_2&=4\pi\rho^2 (w_y^2+ w_x^2)-2\rho( w_y Y_{T0}+ w_x X_{T0}-\rho w_x)\\
    %     a_3&=4\pi^2\rho^2 (w_y^2+ w_x^2)-4\pi\rho(w_y Y_{T0}+ w_x X_{T0}-\rho w_x)+Y_{T0}^2+X_{T0}^2-2\rho X_{T0}
    %     \end{aligned}
    % \right.
    % \nonumber
    % \end{equation}
    
    \item [(2)] The solution path is of $LR_{2\pi}$ if and only if
    \begin{equation}\label{EQ:CC2_1}
    \left\{
        \begin{aligned}
        \theta_f &=\alpha+\pi/2-2n\pi\\
        w_x &= \dfrac {\dfrac{X_{T0}}{\rho}-\cos(\theta_{f} - \dfrac{\pi}{2})+1}{\theta_{f} - \dfrac{\pi}{2}+2(n+1)\pi}\\
        w_y &= \dfrac {\dfrac{Y_{T0}}{\rho}-\sin(\theta_{f} - \dfrac{\pi}{2})}{\theta_{f} - \dfrac{\pi}{2}+2(n+1)\pi}
        \end{aligned}
    \right.
    \end{equation}
    where $n=0$ (resp. $=1$) if $\theta_{f}\geq \pi/2$ (resp. $< \pi/2$); in addition, if Eq.~(\ref{EQ:CC2_1}) holds, we have that $\alpha$ is a root of
    \begin{align}\label{EQ:CC2_2}
     b_1 \alpha^2+b_2\alpha+b_3=0   
    \end{align}
    where the constants $b_1$--$b_3$ are presented in Appendix \ref{Appendix:A}.
    % \begin{equation}
    % \left\{
    %     \begin{aligned}
    %     b_1&=\rho^2 w_y^2+\rho^2 w_x^2\\
    %     b_2&=4\pi\rho^2 (w_y^2+ w_x^2)-2\rho( w_y Y_{T0}+ w_x X_{T0}-\rho w_x)\\
    %     b_3&=4\pi^2\rho^2 (w_y^2+ w_x^2)-4\pi\rho(w_y Y_{T0}+ w_x X_{T0}-\rho w_x)+Y_{T0}^2+X_{T0}^2+2\rho X_{T0}
    %     \end{aligned}
    % \right.
    % \nonumber
    % \end{equation}
    \end{description}
\end{lemma}
The proof of this lemma is postponed to Appendix \ref{Appendix:A}.
As a result of Lemma \ref{Le:CC}, one is able to use Eq.~(\ref{EQ:CC1_1}) and Eq.~(\ref{EQ:CC2_1}) to check if the path of type $CC_{2\pi}$ is feasible or not.
If Eq.~(\ref{EQ:CC1_1}) or Eq.~(\ref{EQ:CC2_1}) is satisfied, one can find $\alpha$ by solving the corresponding analytical equations in Eq.~(\ref{EQ:CC1_2}) or Eq.~(\ref{EQ:CC2_2}). 

Then, by the definition of $T_{CC}$, if the path type is $RL_{2\pi}$, we have
\begin{equation}\label{Eq:T_RL}
T_{RL} = \left\{
\begin{aligned}
&\arg\min\limits_{\alpha \in \{\rm{roots of Eq.~(\ref{EQ:CC1_2})}\}} \left\{ \rho(\alpha+2\pi) \right\}, \quad &&\text{if Eq.~(\ref{EQ:CC1_1}) is satisfied}\\
&+\infty,\quad && \text{if Eq.(\ref{EQ:CC1_1}) is not satisfied}
\end{aligned}
\right.
\end{equation}
In addition, the control along the path of $RL_{2\pi}$ is determined by Eq.~(\ref{Eq:U_CC}) if Eq.(\ref{EQ:CC1_1}) is satisfied. 

Similarly, if the path type is $LR_{2\pi}$, we have
\begin{equation}
T_{LR} = \left\{
\begin{aligned}
&\arg\min\limits_{\alpha \in \{\rm{roots of Eq.~(\ref{EQ:CC2_2})}\}} \left\{ \rho(\alpha+2\pi) \right\}, \quad &&\text{if Eq.~(\ref{EQ:CC2_1}) is satisfied}\\
&+\infty,\quad && \text{if Eq.(\ref{EQ:CC2_1}) is not satisfied}
\end{aligned}
\right.
\end{equation}
In addition, the control along the path of $LR_{2\pi}$ is also determined by Eq.~(\ref{Eq:U_CC}) if Eq.(\ref{EQ:CC2_1}) is satisfied. 

\subsection{Closed-Form Solution for Path of \texorpdfstring{$CCC$}{CCC} }

\begin{figure}[htbp]
    \centering
    \subfigure[$RL_{<\pi}R$]{
        \begin{minipage}[t]{7cm}
        \centering
        \includegraphics[height=4.2cm]{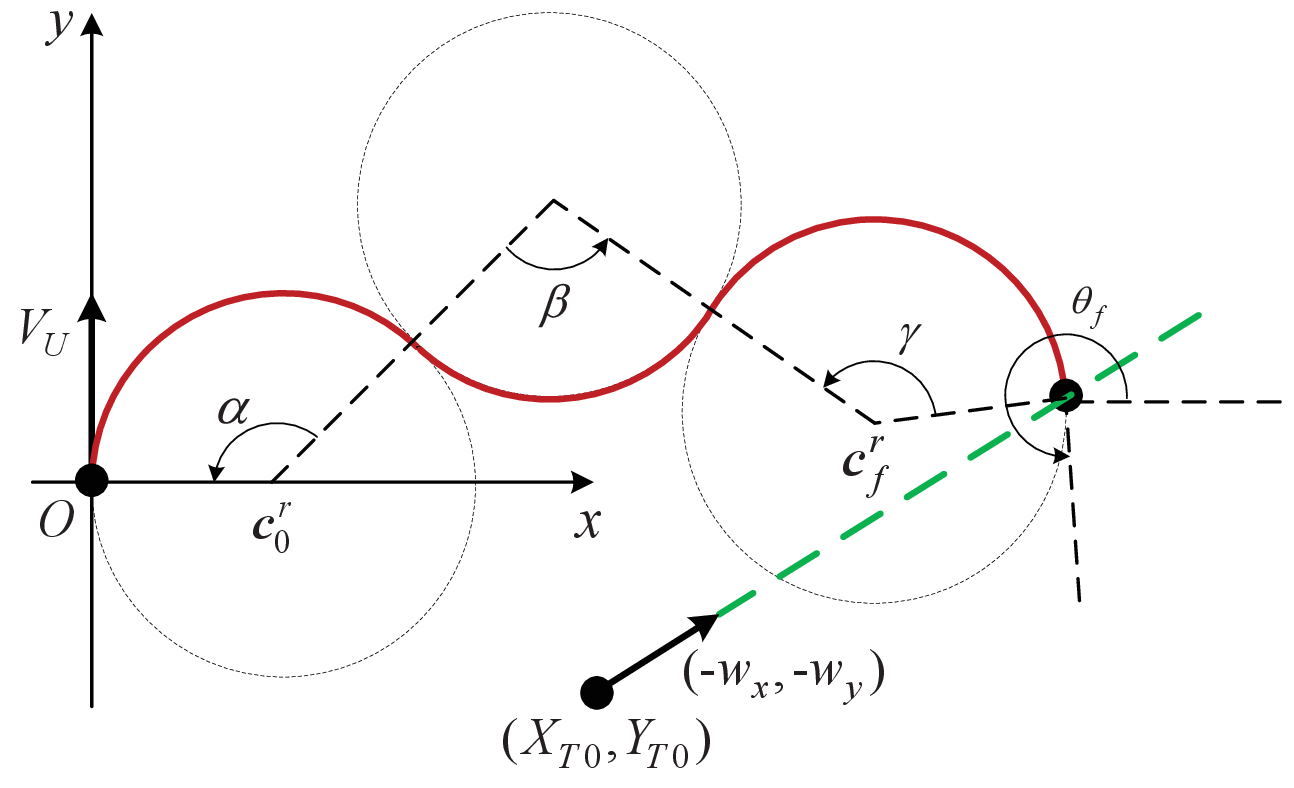}
        \label{Fig:geometry_RLR}
    \end{minipage}%
    }%
    \subfigure[$LR_{<\pi}L2$]{
        \begin{minipage}[t]{7cm}
        \centering
        \includegraphics[height=4.3cm]{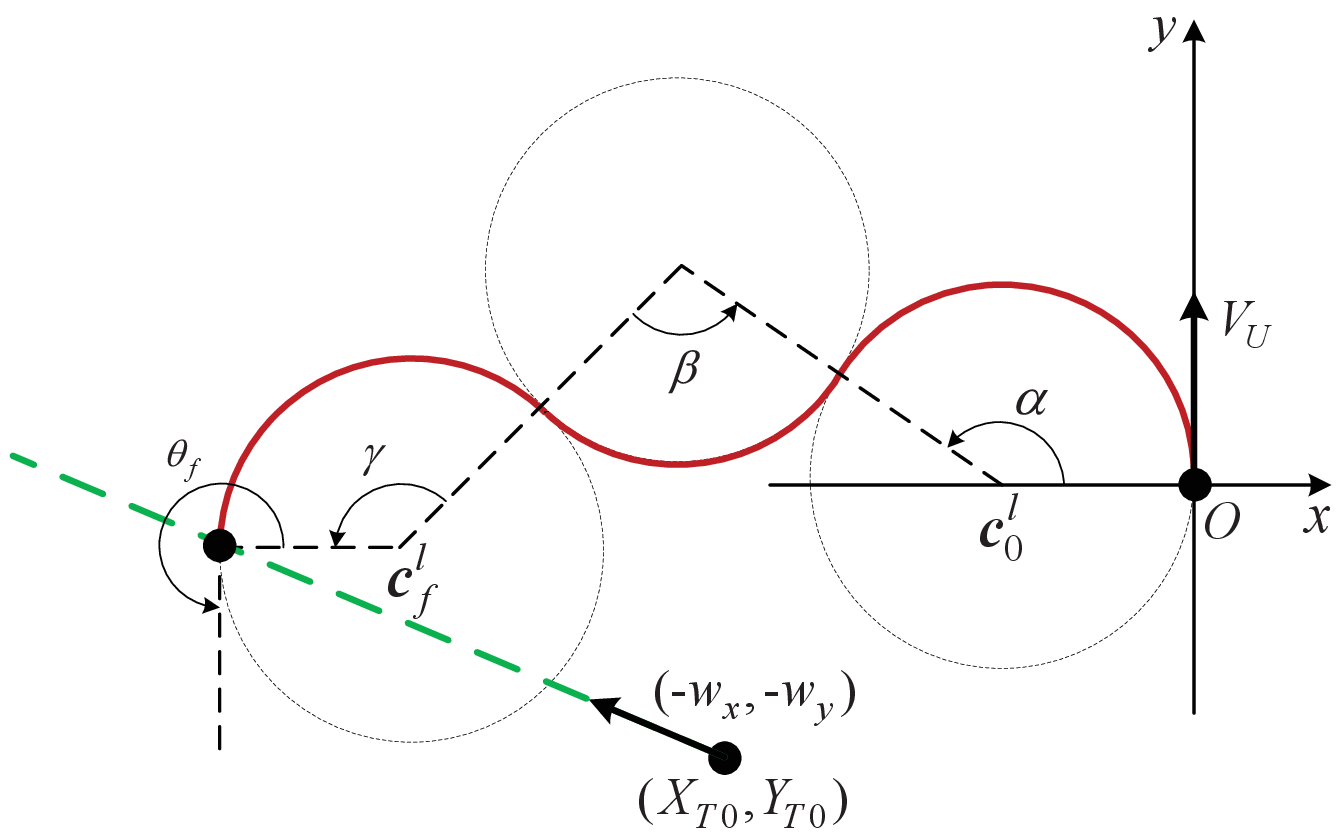}
        \label{Fig:geometry_LRL}
    \end{minipage}%
    }%
    
    \subfigure[$RL_{>\pi}R$]{
        \begin{minipage}[t]{7cm}
        \centering
        \includegraphics[height=4.65cm]{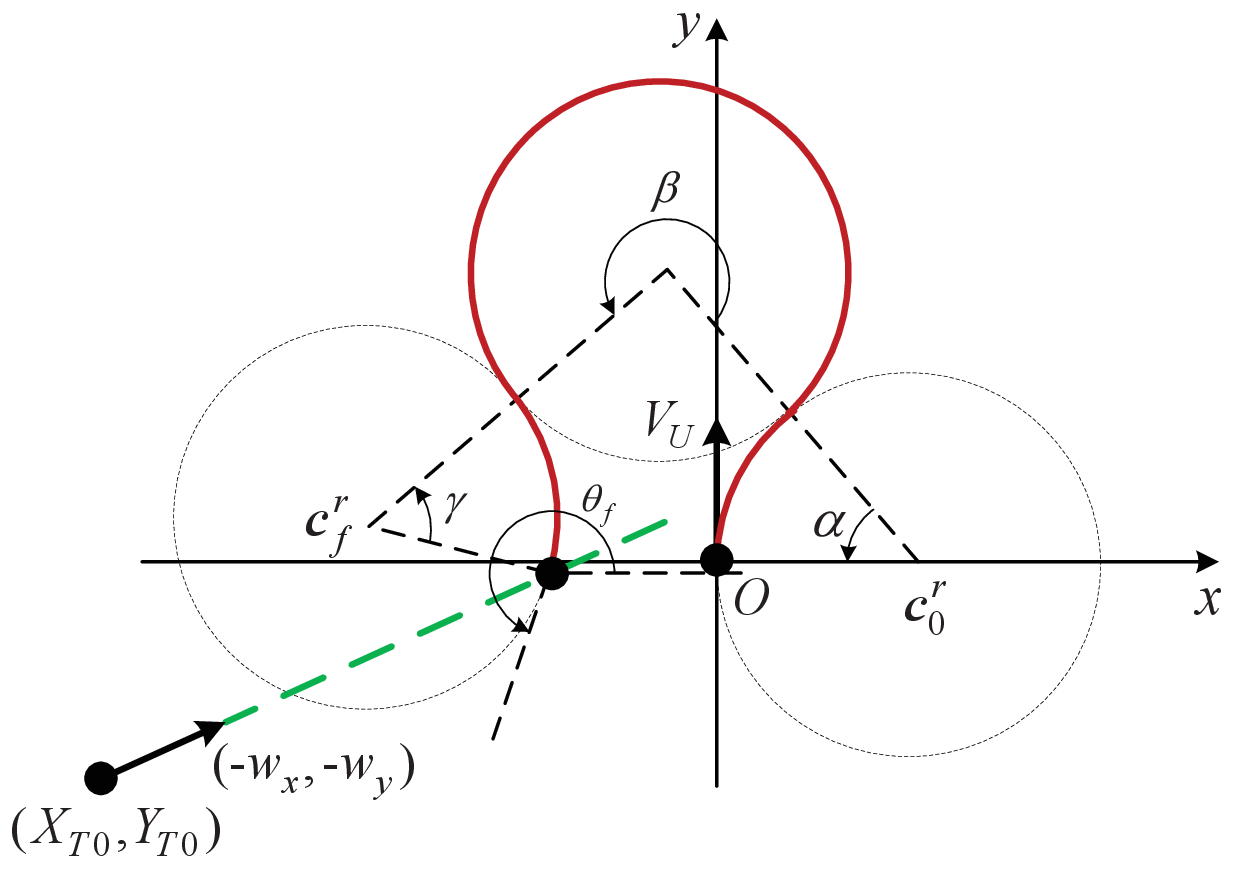}
        \label{Fig:geometry_RLR2}
    \end{minipage}%
    }%
    \subfigure[$LR_{>\pi}L2$]{
        \begin{minipage}[t]{7cm}
        \centering
        \includegraphics[height=4.4cm]{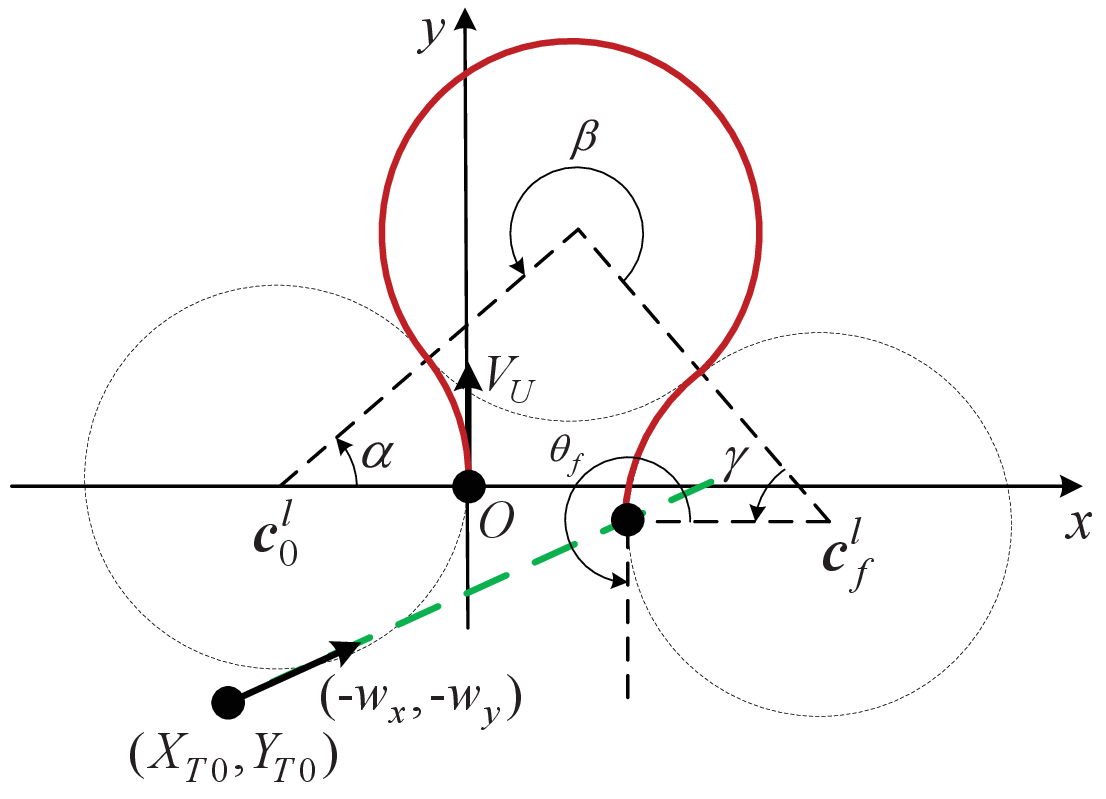}
        \label{Fig:geometry_LRL2}
    \end{minipage}%
    }%
    
    \centering
    \caption{Geometry for the path of $CCC$ in the air-relative frame.}
    \label{Fig:geometryCCC}
\end{figure}

The $CCC$-path includes four different candidates, i.e.,
$$CCC=\{RL_{<\pi}R, RL_{>\pi}R, LR_{<\pi}L, LR_{>\pi}L \}$$ all of which are illustrated in Fig.~\ref{Fig:geometryCCC}.
For notational simplicity, let us denote by $\alpha\in [0,2\pi)$, $\beta\in [0,2\pi)$, and $\gamma\in [0,2\pi)$ the radians of the first, second, and third circular arcs along the $CCC$-path, respectively. Additionally, denote by $U_{CCC}$ the control along the $CCC$-path. Then, we have that $U_{CCC}$ is expressed as
\begin{equation}\label{Eq:U_CCC}
    U_{CCC}(t)=\left\{
    \begin{aligned}
        &\sigma,\quad \quad &&t\in[0,\rho\alpha)\\
        &-\sigma, \quad \quad &&t\in[\rho\alpha,\rho(\alpha+\beta))\\
        &\sigma, \quad \quad &&t\in[\rho(\alpha+\beta),\rho(\alpha+\beta+\gamma)]
    \end{aligned}
    \right.
\end{equation}
where $\sigma=-1$ if the solution path is of $RLR$, and $\sigma=1$ if the solution path is of $LRL$.
It is now enough to find the value of $\alpha$, $\beta$ and $\gamma$ in order to find the solution path of $CCC$.
By the following lemma, we shall show that $\beta$ is a root of a nonlinear equation, and the values of $\alpha$ and $\gamma$ are totally determined by $\beta$.
\begin{lemma}\label{TH}
If the solution path for Problem \ref{problem2} is of type $CCC$, the following two statements hold:
\begin{description}
\item[(1)] if the solution path is of $RLR$, we have
\begin{align}
& c_1\beta^2 + c_2 \beta +c_3\cos \beta +c_4 =0
\label{EQ:RLR}\\
&cos (-\alpha + \dfrac{\beta}{2} + \dfrac{\pi}{2} )= \dfrac{-(X_{T0}+\rho w_x(\dfrac{\pi}{2}- \theta_f+2\beta+2n\pi)) - \rho \sin \theta_f + \rho }{4\rho \sin (-\dfrac{\beta}{2})} 
\label{EQ:RLR_alpha}\\
&cos (\gamma - \dfrac{\beta}{2} - \theta_f) =  \dfrac{-(X_{T0}+\rho w_x(\dfrac{\pi}{2}- \theta_f+2\beta+2n\pi)) - \rho \sin \theta_f + \rho}{4\rho \sin (-\dfrac{\beta}{2})}
\label{EQ:RLR_gamma}
\end{align}
where the positive integer $n$ takes values so that  $(\dfrac{\pi}{2} - \alpha + \beta - \gamma +2n\pi) \in [0,2\pi)$, and the constants $c_1$--$c_4$ in Eq.~(\ref{EQ:RLR}) are presented in Appendix \ref{Appendix:A}.
% \begin{equation}
%     \left\{
%     \begin{aligned}
%     a_1 &=4 \rho^2 (w_x^2+w_y^2)\\
%     a_2 &=4 \rho (M w_x -N w_y)\\
%     a_3 &=8\rho^2\\
%     a_4 &=M^2+N^2-8 \rho^2
%     \end{aligned}\nonumber
%     \right.
%     \label{Eq:a1-a4}
% \end{equation}
% with
% \begin{equation}
%     \left\{
%     \begin{aligned}
%      M &=-X_{T0}+\rho w_x(\dfrac{\pi}{2} - \theta_f + 2n\pi)-\rho \sin \theta_f + \rho \\
%     N &=Y_{T0} - \rho w_y (\dfrac{\pi}{2} - \theta_f + 2n\pi) - \rho \cos \theta_f
%     \end{aligned}\nonumber
%     \right.
%     \label{Eq:M-N}
% \end{equation}
\item[(2)] if the solution path is of $LRL$, we have
\begin{align}
& d_1\beta^2 + d_2 \beta +d_3\cos \beta +d_4 =0
\label{EQ:LRL}\\
&cos(\alpha -\dfrac{\beta}{2} +\dfrac{\pi}{2}) = \dfrac{X_{T0}+\rho w_x(\theta_f+2\beta- \dfrac{\pi}{2}+2n\pi)-\rho \sin \theta_f + \rho }{4 \rho \sin (\dfrac{\beta}{2})}
\label{EQ:LRL_alpha}\\
& cos(-\gamma +\dfrac{\beta}{2} +\theta_f)=  \dfrac{X_{T0}+\rho w_x(\theta_f+2\beta- \dfrac{\pi}{2}+2n\pi)-\rho \sin \theta_f + \rho}{4 \rho \sin (\dfrac{\beta}{2})}
 \label{EQ:LRL_gamma}
\end{align}
where the positive integer $n$ takes values so that $(\dfrac{\pi}{2} + \alpha - \beta + \gamma +2n\pi) \in [0,2\pi)$, and the constants $d_1$--$d_4$ in Eq.~(\ref{EQ:LRL}) are presented in Appendix \ref{Appendix:A}.
% \begin{equation}
%     \left\{
%     \begin{aligned}
%     b_1 &=4 \rho^2 (w_x^2+w_y^2)\\
%     b_2 &=4 \rho (P w_x-Q w_y)\\
%     b_3 &=8\rho^2\\
%     b_4 &=P^2+Q^2-8 \rho^2
%     \end{aligned}\nonumber
%     \right.
%     \label{Eq:b1-b4}
% \end{equation}
% with
% \begin{equation}
%     \left\{
%     \begin{aligned}
%     P &=X_{T0}-\rho w_x(\theta_f - \dfrac{\pi}{2} + 2n\pi)-\rho \sin \theta_{f} + \rho\\
%     Q &= - Y_{T0} +\rho w_y (\theta_f - \dfrac{\pi}{2} + 2n\pi) - \rho \cos \theta_f
%     \end{aligned}\nonumber
%     \right.
% \end{equation}
\end{description}
\end{lemma}
The proof of this lemma is postponed to Appendix \ref{Appendix:A}. 

Lemma \ref{TH} indicates that, to find the value of  $\beta$, it involves solving the transcendental equations in Eq.~(\ref{EQ:RLR}) and Eq.~(\ref{EQ:LRL}). These equations may have multiple roots, but only one specific root is relevant to the optimal path. Existing numerical methods, such as the Newton-Raphson method and the bisection method, cannot guarantee to find
the desired root related to the optimal path. In Appendix \ref{Appendix:B}, we introduce an improved bisection method so that all the roots of Eq.~(\ref{EQ:RLR}) and Eq.~(\ref{EQ:LRL}) can be identified robustly and efficiently. Note that the values of $\alpha$ and $\gamma$ corresponding to each real root $\beta$ can be computed through Eqs.~(\ref{EQ:RLR_alpha}--\ref{EQ:RLR_gamma}) or Eqs.~(\ref{EQ:LRL_alpha}--\ref{EQ:LRL_gamma}). 

Then, by the definition of $T_{CCC}$, if the path type is $RLR$, we have
\begin{equation}
T_{RLR} = \left\{
\begin{aligned}
&\rho(\alpha+\beta+\gamma), \quad &&\text{if the root for Eq.~(\ref{EQ:RLR}) exists}\\
&+\infty,\quad && \text{if the root for Eq.~(\ref{EQ:RLR}) does not exist}
\end{aligned}
\right.
\end{equation}
where $\beta$ is the root for Eq.~(\ref{EQ:RLR}), and $\alpha$ and $\gamma$ are computed through Eqs.~(\ref{EQ:RLR_alpha}--\ref{EQ:RLR_gamma}). In addition, the control along the path of $RLR$ is determined by Eq.~(\ref{Eq:U_CCC}) if the root for Eq.~(\ref{EQ:RLR}) exists. 

Similarly, if the path type is $LRL$, we have
\begin{equation}
T_{LRL} = \left\{
\begin{aligned}
&\rho(\alpha+\beta+\gamma) , \quad &&\text{if the root for Eq.~(\ref{EQ:LRL}) exists}\\
&+\infty,\quad && \text{if the root for Eq.~(\ref{EQ:LRL}) does not exist}
\end{aligned}
\right.
\end{equation}
where $\beta$ is the root for Eq.~(\ref{EQ:LRL}), and $\alpha$ and $\gamma$ are computed through Eqs.~(\ref{EQ:LRL_alpha}--\ref{EQ:LRL_gamma}). The control along the path of $LRL$ is also determined by Eq.~(\ref{Eq:U_CCC}) if the root for Eq.~(\ref{EQ:LRL}) exists. 

\subsection{Closed-Form Solution for Path of \texorpdfstring{$CSC$}{CSC} }

The $CSC$-path includes four different candidates, i.e.,
$$CSC=\{RSR, RSL, LSR, LSL \}$$
For notational simplicity, let us denote by $\alpha\in [0,2\pi)$ and $\gamma\in [0,2\pi)$ the radians of the first and the third circular arcs along the $CSC$-path, respectively.
Additionally, denote by $d$ and $\beta$ the length and orientation angle of the straight-line segment along the $CSC$-path.
Denote by $U_{CSC}$ the control along the $CSC$-path. Then, we have that $U_{CSC}$ is expressed as
\begin{equation}\label{Eq:U_CSC}
    U_{CSC}(t)=\left\{
    \begin{aligned}
        &\sigma,\quad \quad &&t\in[0,\rho\alpha)\\
        &0, \quad \quad &&t\in[\rho\alpha,\rho\alpha+d)\\
        &\kappa, \quad \quad &&t\in[\rho\alpha+d,\rho(\alpha+\gamma)+d]
    \end{aligned}
    \right.
\end{equation}
where $\sigma=-1$ and $\kappa=-1$ for the $RSR$-path, $\sigma=-1$ and $\kappa=1$ for the $RSL$-path, $\sigma=1$ and $\kappa=-1$ for the $LSR$-path, and $\sigma=1$ and $\kappa=1$ for the $LSL$-path.
It is now enough to
find the value of $\alpha$, $d$, and $\gamma$ in order to find the solution path of $CSC$.

It has been shown in \cite{ZHENG2021} that the values of $\alpha$, $d$ and $\gamma$ are entirely determined by $\beta$ for the $CSC$-path, and the closed-form expressions are given as well. For the completeness of the current paper, the closed-form solutions for $\beta$ are summarized in the following lemma. 
\begin{lemma}[Zheng {\it et al.}~\cite{ZHENG2021}]
\label{LE:Zheng}
If the solution path for Problem \ref{problem2} is of type $CSC$, the following two statements hold:
	\begin{description}
		\item[(1)] If the solution path is of $RSR$ or $LSL$, we have
		\begin{align}\label{EQ:RSR}
			e_1 + e_2 \sin \beta +e_3\cos \beta  =0
		\end{align}
  where $e_1$--$e_3$ are presented in Appendix \ref{Appendix:A}.
		\item[(2)] If the solution path is of $RSL$ or $LSR$, we have
		\begin{align}\label{EQ:RSL}
			f_1 + f_2 \sin \beta +f_3 \cos \beta + \beta(f_4 \sin \beta +f_5 \cos \beta) = 0
		\end{align}
  where $f_1$--$f_5$ are presented in Appendix \ref{Appendix:A}.
	\end{description}
\end{lemma}

In \cite{ZHENG2021}, Zheng {\it et al.} proposed a method for finding all the roots of Eq.~(\ref{EQ:RSR}) and Eq.~(\ref{EQ:RSL}). Then, the value of $\alpha$, $d$, and $\gamma$ corresponding to each real root $\beta$ can be computed analytically, and readers who are interested in those analytical expressions are referred to \cite{ZHENG2021}.

Then, by the definition of $T_{CSC}$, if the path type is $RSR$ or $LSL$, we have
\begin{equation}
T_{RSR \text{ or } LSL} = \left\{
\begin{aligned}
&\rho(\alpha+\gamma)+d, \quad &&\text{if the root for Eq.~(\ref{EQ:RSR}) exists}\\
&+\infty,\quad && \text{if the root for Eq.~(\ref{EQ:RSR}) does not exist}
\end{aligned}
\right.
\end{equation}
In addition, the control along the path of $RSR$ or $LSL$ is determined by Eq.~(\ref{Eq:U_CSC}) if the root for Eq.~(\ref{EQ:RSR}) exists. 

Similarly, if the path type is $RSL$ or $LSR$, we have
\begin{equation}
T_{RSL \text{ or } LSR} = \left\{
\begin{aligned}
&\rho(\alpha+\gamma)+d, \quad &&\text{if the root for Eq.~(\ref{EQ:RSL}) exists}\\
&+\infty,\quad && \text{if the root for Eq.~(\ref{EQ:RSL}) does not exist}
\end{aligned}
\right.
\end{equation}
and the control along the path of $RSL$ or $LSR$ is also determined by Eq.~(\ref{Eq:U_CSC}) if the root for Eq.~(\ref{EQ:RSL}) exists.

\section{Numerical Simulations}\label{SE:simulation}

In this section, we shall present two numerical examples to illustrate how to use the closed-form solutions derived in Section \ref{SE:Analytical} to find the globally optimal solution in real time. 

\subsection{Case 1}

\begin{figure}[htbp]
    \centering
    \subfigure[Paths in air-relative frame]{
    \begin{minipage}[t]{\textwidth}
        \centering
        \includegraphics[width=.65\textwidth]{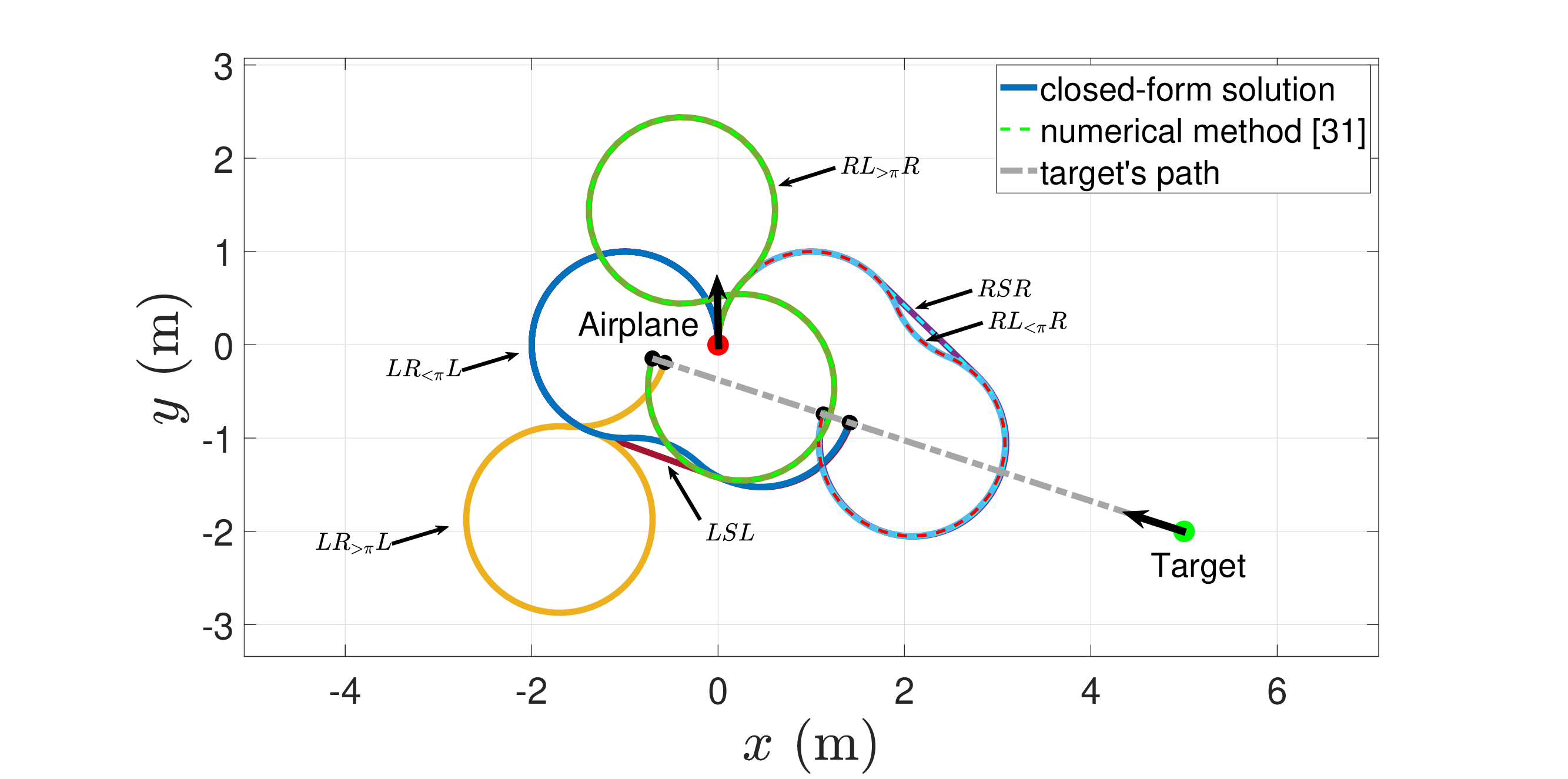}
        \label{Fig:case1_tra1}
    \end{minipage}
    }%
    \newline
    \subfigure[Paths in frame $OXY$]{
    \begin{minipage}[t]{\textwidth}
        \centering
        \includegraphics[width=.65\textwidth]{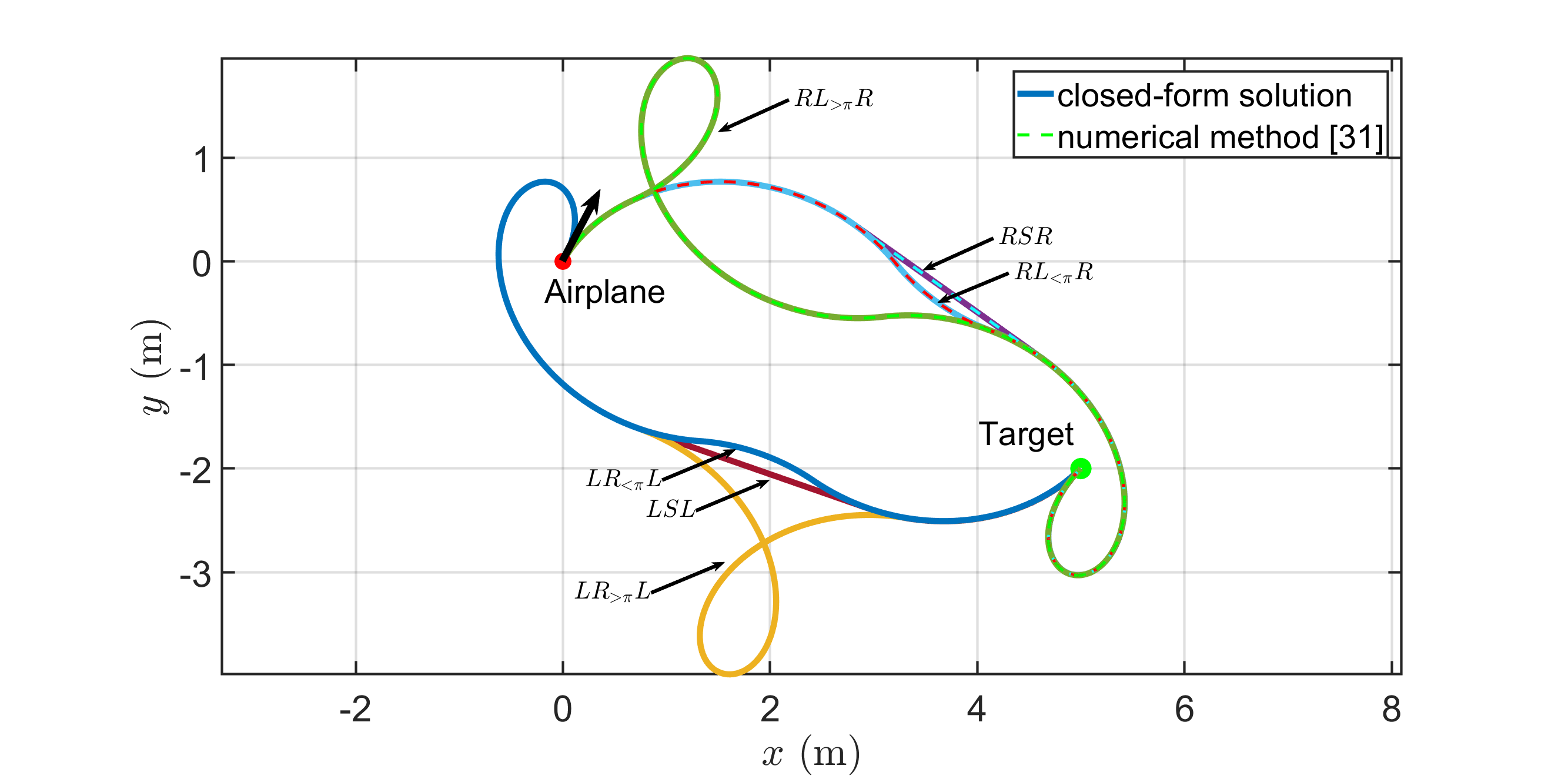}
        \label{Fig:case1_tra2}
    \end{minipage}
    }%
    \caption{Paths of all candidate types for Case 1.}
    \label{Fig:type4}
\end{figure}

Consider a Dubins airplane moving in a steady wind, where the velocity of the wind is $(w_x,w_y)=(0.475,-0.155)$ m/s. The initial state of the Dubins airplane is $(0,0,\pi/2)$, and the speed of the airplane is 1 m/s. 
The target is located at $(X_{T0},Y_{T0})=(5,-2)$ m. Set the parameter $\rho$ as 1. The Dubins airplane is expected to arrive the target with a terminal heading angle of 72 deg while minimizing the time duration $t_f$. The paths for all candidate types generated by the closed-form solutions in Section \ref{SE:Analytical} and those by the numerical method in \cite{BUZIKOV2022109968} are presented by solid curves and dotted curves in Fig.~\ref{Fig:type4}, respectively. 

\begin{figure}[hbt!]
    \centering
    \includegraphics[width=.65\textwidth]{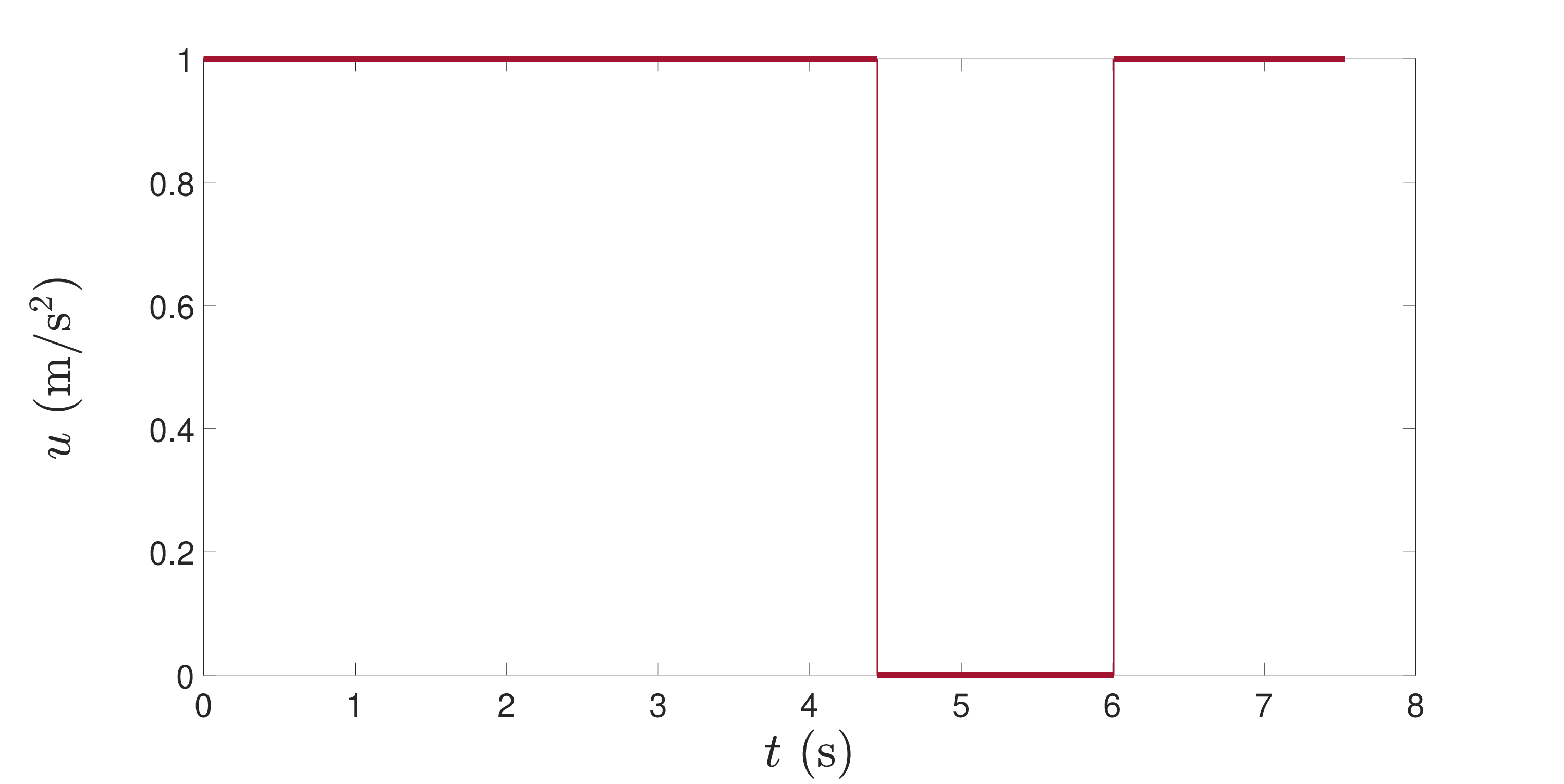}
    \caption{Optimal Control Profile for Case 1.}
    \label{Fig:type4_control}
\end{figure}

From Fig.~\ref{Fig:type4}, we can see that 6 candidate paths (of $RL_{>\pi}R$, $RL_{<\pi}R$, $LR_{>\pi}L$, $LR_{<\pi}L$, $RSR$, $LSL$) are obtained by the closed-form solution. 
The time durations for the paths of $RL_{>\pi}R$, $RL_{<\pi}R$, $LR_{>\pi}L$, $LR_{<\pi}L$, $RSR$, $LSL$ are 11.9937 s, 8.1420 s, 11.7152 s, 7.5570 s, 8.1157 s, 7.5294 s, respectively. It is evident that the time duration of the $LSL$-path is the smallest. Therefore, according to Eq.~(\ref{Eq:t_f}), we have that the global minimum time for reaching the target is 7.5294 s, and the type of globally optimal path is of $LSL$. By applying Eq.~(\ref{Eq:U_CSC}), the corresponding control profile of the $LSL$-path is presented in Fig.~\ref{Fig:type4_control}.

It is worth noting that only 3 candidate paths are found by the numerical method in \cite{BUZIKOV2022109968}, and the $LSL$-path does not belong to the three. This indicates that the numerical method in \cite{BUZIKOV2022109968} cannot guarantee to find the globally optmal solution. Compared to the numerical methods in \cite{BUZIKOV2022109968}, the closed-form solution proposed in this paper guarantees finding the globally minimum-time path for Problem \ref{problem2}.

\subsection{Case 2}

\begin{figure}[htbp]
    \centering
    \subfigure[Paths in the air-relative frame]{
    \begin{minipage}[t]{\textwidth}
        \centering
        \includegraphics[width=.65\textwidth]{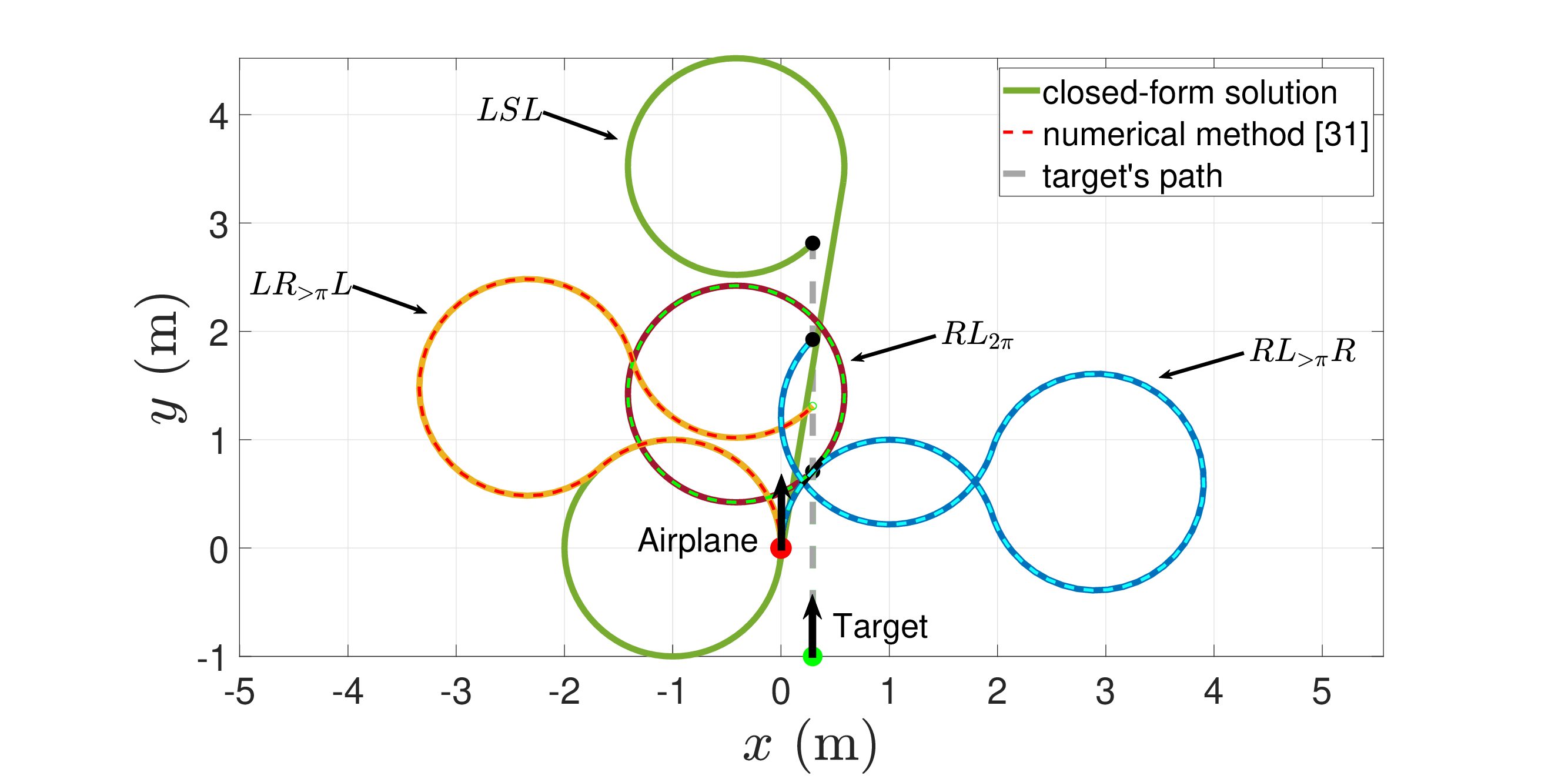}
        \label{Fig:case2_tra1}
    \end{minipage}
    }%
    \newline
    \subfigure[Paths in frame $OXY$]{
    \begin{minipage}[t]{\textwidth}
        \centering
        \includegraphics[width=.65\textwidth]{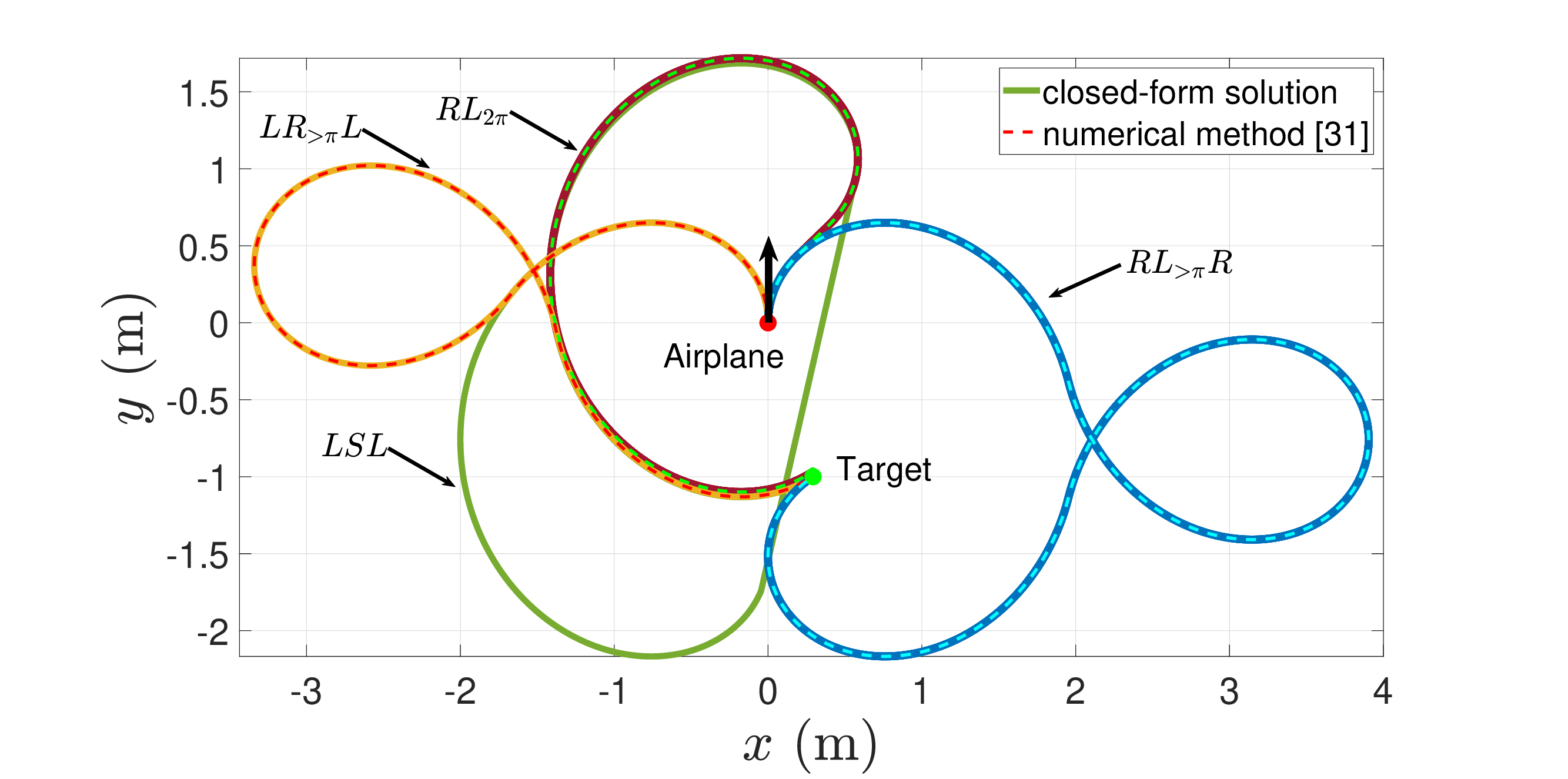}
        \label{Fig:case2_tra2}
    \end{minipage}
    }%
    \caption{Paths of all candidate types for Case 2.}
    \label{Fig:type7}
\end{figure}

Consider a Dubins airplane moving in a steady wind,  with initial conditions identical to the case in Fig.~9 of \cite{BUZIKOV2022109968}. The velocity of the wind is $(w_x,w_y)=(0,-\dfrac{4+2\sqrt{2}}{9\pi})$ m/s. The initial state of the Dubins airplane is $(0,0,\pi/2)$, and the speed of the airplane is 1 m/s. 
The target is located at $(X_{T0},Y_{T0})=(1-\dfrac{1}{\sqrt{2}}, -1)$ m. Set the parameter $\rho$ as 1. The Dubins airplane is expected to arrive the target with a terminal heading angle of 45 deg while minimizing the time duration $t_f$. The paths for all candidate types generated by the closed-form solutions in Section \ref{SE:Analytical} and those by the numerical method in \cite{BUZIKOV2022109968} are presented by solid curves and dotted curves in Fig.~\ref{Fig:type7}, respectively. 

\begin{figure}[hbt!]
    \centering
    \includegraphics[width=.65\textwidth]{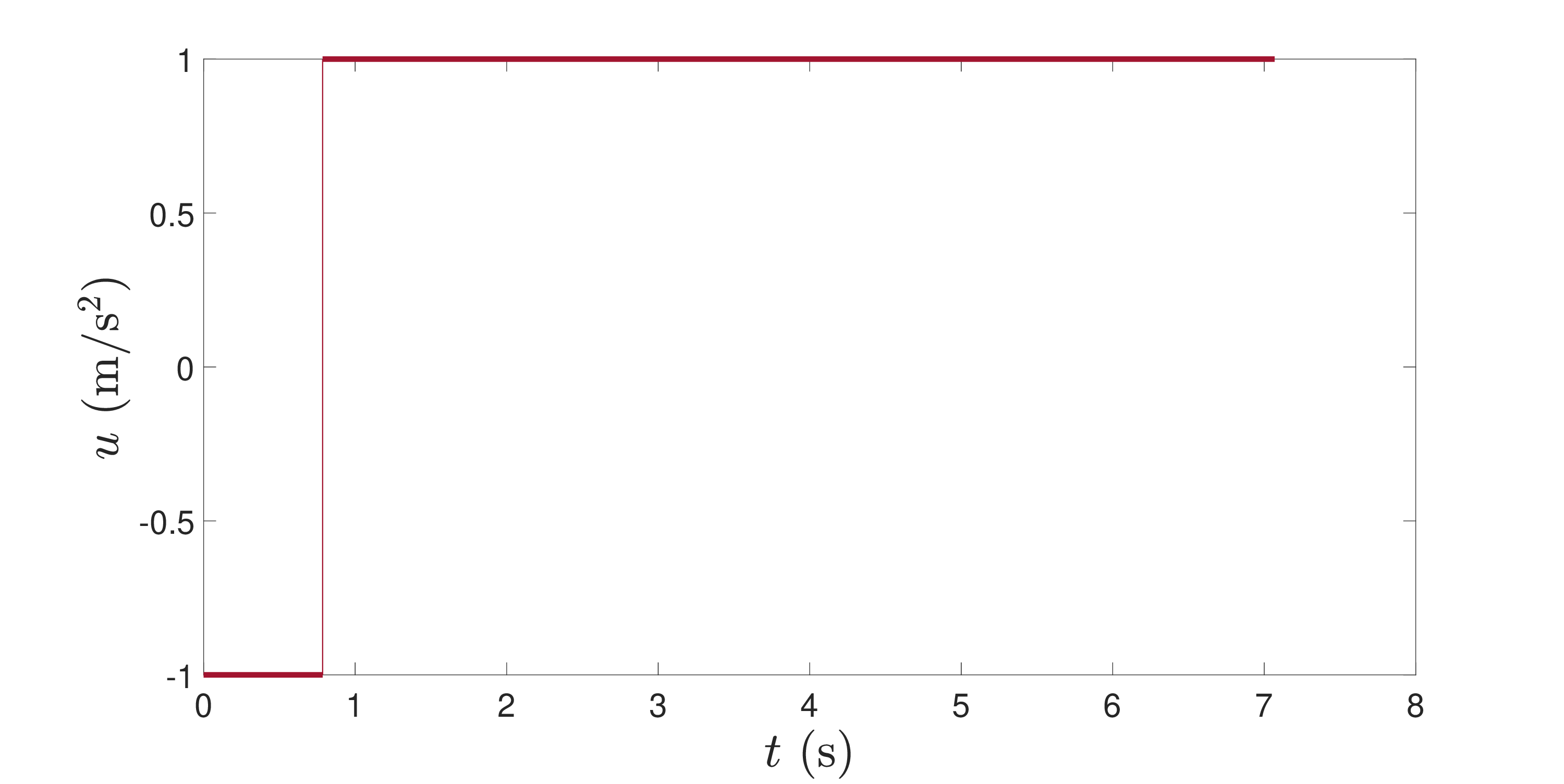}
    \caption{Optimal Control Profile for Case 2.}
    \label{Fig:type7_control}
\end{figure}

From Fig.~\ref{Fig:type7}, we can see that 4 candidate paths (of $LSL$, $LR_{>\pi}L$, $RL_{>\pi}R$, $RL_{2\pi}$) are obtained by the closed-form solution, while only 3 candidate paths are found by the numerical method in \cite{BUZIKOV2022109968}. The time durations for the $LSL$, $LR_{>\pi}L$, $RL_{>\pi}R$ and $RL_{2\pi}$-paths are 15.7929 s, 9.5686 s, 12.1137 s and 7.0686 s, respectively. Thus, according to Eq.~(\ref{Eq:t_f}), we have that the global minimum time is 7.0686 s, and the type of the globally minimum-time path is of $RL_{2\pi}$. 
By applying Eq.~(\ref{Eq:U_CC}), the corresponding control profile of the $RL_{2\pi}$-path is presented in Fig.~\ref{Fig:type7_control}.

\section{Conclusions}\label{SE:conclusion}

This paper addressed the minimum-time path panning problem for a Dubins airplane from an initial configuration to a terminal configuration under the influence of steady wind. First, the minimum-time path-planning problem was formulated as an optimal control problem. Then, necessary conditions derived from the PMP showed that the minimum-time solution for the optimal control problem belongs to a family of four types. The geometric properties were analysed for each type of paths. Specific initial constraints were provided in the paper for the existence of $SC_{2\pi}$- and $CC_{2\pi}$-paths, and explicit expressions of control and path were provided. For the remaining two types $CCC$ and $CSC$, it was found that the paths are determined by some transcendental equations. An improved bisection method was presented to find all roots of those transcendental equations. As a result, the paths of $CCC$ and $CSC$ can be obtained within a constant time. 
Finally, it was demonstrated through numerical simulations that the globally optimal solution can be efficiently obtained by comparing all candidate paths, while the existing numerical methods in the literature fails to find the globally optimal solution.

\section*{Acknowledgments}

This work was supported by the National Natural Science Foundation of China (Nos. 61903331 and 62088101).

\bibliography{sample}

\appendix

\section{Proofs and Supplementary Notes for Lemmas in Section \ref{SE:Analytical}}
\label{Appendix:A}

\subsection*{A.1 Proof of Lemma \ref{Le:CC}}

Let us start by proving the first statement.
According to Fig.~\ref{Fig:geometry_RL2pi}, we have that the solution path is of $RL_{2\pi}$ if and only if 
\begin{equation}
\left\{
    \begin{aligned}
    &x_f=\rho+\rho\cos(\pi-\alpha)\\
    &y_f=\rho\sin(\pi-\alpha)\\
    &\theta_f=\pi/2-\alpha+2n\pi\\
    &-w_x(x_f-X_{T0}) =-w_y(y_f-Y_{T0})
    \end{aligned}
\right.
\label{EQ:CC1}
\end{equation}
where $n=0$ (resp. $=1$) if $\theta_{f}<\pi/2$ (resp. $\geq \pi/2$). 
Note that the time required for the Dubins airplane to travel from its initial state $(x_0,y_0,\theta_0)$ to the desired terminal position $(x_f,y_f)$ is equal to the time required for the target to reach the same terminal position in the air-relative frame, indicating
\begin{align}
	 \rho(\alpha +2\pi)= \dfrac{\sqrt{(x_f - X_{T0})^2 + (y_f - Y_{T0})^2}}{\sqrt{(-w_x)^2 + (-w_y)^2}}
	\label{EQ:CC3}
\end{align}
Combining Eq.~(\ref{EQ:CC1}) and Eq.~(\ref{EQ:CC3}), we have
\begin{equation}
	\left\{
	    \begin{aligned}
		-\rho (\alpha+2\pi)w_x &=x_f-X_{T0}\\
		-\rho (\alpha+2\pi)w_y &=y_f-Y_{T0}
		\end{aligned}
	\right.
	\label{EQ:CC4}
\end{equation}
Squaring both sides of each equation in Eq.~(\ref{EQ:CC4}) and adding the results yields
\begin{equation}
    \label{EQ:CC5}
    \left[Y_{T0}-\rho(\alpha+2\pi)w_y\right]^2+\left[X_{T0}-\rho-\rho(\alpha+2\pi)w_x\right]^2=\rho^2
\end{equation}
Rewriting Eq.~(\ref{EQ:CC5}) leads to
\begin{equation}\label{EQ:CC6}
    a_1\alpha^2+a_2\alpha+a_3=0
\end{equation}
where
\begin{equation}
    \left\{
        \begin{aligned}
        a_1&=\rho^2 w_y^2+\rho^2 w_x^2\\
        a_2&=4\pi\rho^2 (w_y^2+ w_x^2)-2\rho( w_y Y_{T0}+ w_x X_{T0}-\rho w_x)\\
        a_3&=4\pi^2\rho^2 (w_y^2+ w_x^2)-4\pi\rho(w_y Y_{T0}+ w_x X_{T0}-\rho w_x)+Y_{T0}^2+X_{T0}^2-2\rho X_{T0}
        \end{aligned}
    \right.
    \nonumber
    \end{equation}
Solving Eq.~(\ref{EQ:CC6}) to find the value of $\alpha$ and substituting it into Eq.~(\ref{EQ:CC1}) to eliminate $x_f$ and $y_f$, we finally obtain the initial conditions that the $RL_{2\pi}$-path must be satisfied, as shown in Eq.~(\ref{EQ:CC1_1}).

%%%%%%%%%%%%%%%%%%%%%%%%%%%%%%%%%%%%%%%%%%%%%%%%%%%%%%
Let us now proceed to prove the second statement.
According to Fig.~\ref{Fig:geometry_LR2pi}, we have that the solution path is of $LR_{2\pi}$ if and only if 
\begin{equation}
\left\{
    \begin{aligned}
    &x_f=-\rho+\rho\cos\alpha\\
    &y_f=\rho\sin\alpha\\
    &\theta_f=\alpha+\pi/2-2n\pi\\
    &-w_x(x_f-X_{T0}) =-w_y(y_f-Y_{T0})
    \end{aligned}
\right.
\label{EQ:CC1_2_2}
\end{equation}
where $n=0$ (resp. $=1$) if $\theta_{f}\geq\pi/2$ (resp. $< \pi/2$).
Note that the time required for the Dubins airplane to travel from its initial state $(x_0,y_0,\theta_0)$ to the desired terminal position $(x_f,y_f)$ is equal to the time required for the target to reach the same terminal position in the air-relative frame, indicating
\begin{align}
	 \rho(\alpha +2\pi)= \dfrac{\sqrt{(x_f - X_{T0})^2 + (y_f - Y_{T0})^2}}{\sqrt{(-w_x)^2 + (-w_y)^2}}
	\label{EQ:CC3_2}
\end{align}
Combining Eq.~(\ref{EQ:CC1_2_2}) and Eq.~(\ref{EQ:CC3_2}), we have
\begin{equation}
	\left\{
	    \begin{aligned}
		-\rho (\alpha+2\pi)w_x &=x_f-X_{T0}\\
		-\rho (\alpha+2\pi)w_y &=y_f-Y_{T0}
		\end{aligned}
	\right.
	\label{EQ:CC4_2}
\end{equation}
Squaring both sides of each equation in Eq.~(\ref{EQ:CC4_2}) and adding the results yields
\begin{equation}
    \label{EQ:CC5_2}
    \left[Y_{T0}-\rho(\alpha+2\pi)w_y\right]^2+\left[X_{T0}+\rho-\rho(\alpha+2\pi)w_x\right]^2=\rho^2
\end{equation}
Rewriting Eq.~(\ref{EQ:CC5_2}) leads to
\begin{equation}\label{EQ:CC6_2}
    b_1\alpha^2+b_2\alpha+b_3=0
\end{equation}
where
\begin{equation}
    \left\{
        \begin{aligned}
        b_1&=\rho^2 w_y^2+\rho^2 w_x^2\\
        b_2&=4\pi\rho^2 (w_y^2+ w_x^2)-2\rho( w_y Y_{T0}+ w_x X_{T0}-\rho w_x)\\
        b_3&=4\pi^2\rho^2 (w_y^2+ w_x^2)-4\pi\rho(w_y Y_{T0}+ w_x X_{T0}-\rho w_x)+Y_{T0}^2+X_{T0}^2+2\rho X_{T0}
        \end{aligned}
    \right.
    \nonumber
    \end{equation}
Solving Eq.~(\ref{EQ:CC6_2}) to find the value of $\alpha$ and substituting it into Eq.~(\ref{EQ:CC1_2}) to eliminate $x_f$ and $y_f$, we finally obtain the initial conditions that the $LR_{2\pi}$-path must be satisfied, as shown in Eq.~(\ref{EQ:CC2_1}). This completes the proof.

%%%%%%%%%%%%%%%%%%%%%%%%%%%%%%%%%%%%%%%%%%%%%%%%%%%%%%%%%%
\subsection*{A.2 Proof of Lemma \ref{TH}}
Let us start by proving the first statement.
The time required for the Dubins airplane to travel from its initial state $(x_0,y_0,\theta_0)$ to the desired terminal position $(x_f,y_f)$ is equal to the time required for the target to reach the same terminal position in the air-relative frame, indicating
\begin{align}
	\rho (\alpha + \beta + \gamma) =\dfrac{\sqrt{(x_f -X_{T0})^2 + (y_f - Y_{T0})^2}}{\sqrt{(-w_x)^2 + (-w_y)^2}}
	\label{Eq:RLR0}
\end{align}
According to Fig.~\ref{Fig:geometry_RLR}, we have 
\begin{equation}
   \left\{
    \begin{aligned}
    &\dfrac{\pi}{2} - \alpha + \beta - \gamma +2n\pi= \theta_f \\
    &\boldsymbol{c}_0^r + 2\rho 
    \left[
    \begin{array}{c}
    \cos(\pi -\alpha)\\
    \sin(\pi -\alpha)
    \end{array}
    \right]=
    \boldsymbol{c}_f^r + 2\rho 
    \left[
    \begin{array}{c}
    \cos(\theta_f + \gamma + \dfrac{\pi}{2})\\
    \sin(\theta_f + \gamma + \dfrac{\pi}{2})
    \end{array}
    \right]
    \end{aligned}
    \label{EQ:theta0_thetaf}
    \right.
\end{equation}
where the positive integer $n$ takes values so that  $(\dfrac{\pi}{2} - \alpha + \beta - \gamma +2n\pi) \in [0,2\pi)$, and $\boldsymbol{c}_0^r$ and $\boldsymbol{c}_f^r$ are the centers of the first arc and the third arc, respectively.
Rearranging the second equation in Eq.~(\ref{EQ:theta0_thetaf}) yields
\begin{equation}
	\left\{
	\begin{aligned}
		&2\rho [\sin (\dfrac{\pi}{2} - \alpha) - \sin(\theta_f + \gamma)] = - x_f - \rho \sin \theta_{f} + \rho \\
        &2\rho [\cos (\dfrac{\pi}{2} - \alpha) - \cos (\theta_f + \gamma)]= y_f - \rho \cos \theta_{f} 
	\end{aligned}
	\label{Eq:RLR5}
	\right.
\end{equation} 
Combining the first equation in Eq.~(\ref{EQ:theta0_thetaf}) and Eq.~(\ref{Eq:RLR5}), we have
\begin{equation}
	\left\{
	\begin{aligned}
		& 4\rho \cos \left(\dfrac{\pi/2 - \alpha  +\theta_f + \gamma}{2}\right) \sin (-\dfrac{\beta}{2})  = - x_f - \rho \sin \theta_f + \rho \\
   & -4\rho \sin \left(\dfrac{\pi/2 - \alpha  +\theta_f + \gamma}{2}\right) \sin (-\dfrac{\beta}{2})  = y_f  - \rho \cos \theta_{ f} 
	\end{aligned}
	\label{Eq:RLR6}
	\right.
\end{equation} 
Squaring both sides of each equation in Eq.~(\ref{Eq:RLR6}) and adding the results yields
\begin{align}
16 \rho^2 \sin^2 \dfrac{\beta}{2} = \left[ - x_f - \rho \sin \theta_f + \rho \right]^2 + \left[ y_f  - \rho \cos \theta_f \right]^2
\label{Eq:RLR7}
\end{align}
Substituting Eq.~(\ref{EQ:theta0_thetaf}) into Eq.~(\ref{Eq:RLR0}), we have
\begin{align}
\rho(\dfrac{\pi}{2} - \theta_f + 2n\pi) + 2 \rho \beta = \dfrac{\sqrt{(x_f -X_{T0})^2 + (y_f -Y_{T0})^2}}{\sqrt{w_x^2 + w_y^2}}
\label{Eq:RLR8}
\end{align}
Notice that
\begin{align}
    -w_y(x_f-X_{T0}) =-w_x(y_f-Y_{T0})
    \label{Eq:RLR4}
\end{align}
Combining Eq.~(\ref{Eq:RLR4}) and Eq.~(\ref{Eq:RLR8}) leads to
\begin{equation}
\left\{
\begin{aligned}
x_f=X_{T0}+\rho w_x(\dfrac{\pi}{2}- \theta_{ f}+2\beta+2n\pi)\\
y_f = Y_{T0} +\rho w_y(\dfrac{\pi}{2} - \theta_f + 2\beta + 2n \pi)
    \label{Eq:RLR9}
\end{aligned}
\right.
\end{equation}
Substituting Eq.~(\ref{Eq:RLR4}) and Eq.~(\ref{Eq:RLR9}) into Eq.~(\ref{Eq:RLR7}) to eliminate $x_f$ and $y_f$, we obtain a transcendental equation, i.e.,
\begin{align}
    16\rho^2\sin^2\dfrac{\beta}{2}=(M+2\rho w_x\beta)^2+(N-2\rho w_y\beta)^2
\label{Eq:RLR10}
\end{align}
where
\begin{equation}
    \left\{
    \begin{aligned}
    M &=-X_{T0}-\rho w_x(\dfrac{\pi}{2} - \theta_{f} + 2n\pi)-\rho \sin \theta_f + \rho \\
    N &=Y_{T0} + \rho w_y (\dfrac{\pi}{2} - \theta_{f} + 2n\pi) - \rho \cos \theta_f 
    \end{aligned}\nonumber
    \right.
\end{equation}
Rearranging Eq.~(\ref{Eq:RLR10}), we have 
\begin{align}
    c_1\beta^2 + c_2 \beta +c_3\cos \beta +c_4 =0
\end{align}
where $c_1$--$c_4$ are given as
\begin{equation}
    \left\{
    \begin{aligned}
    c_1 &=4 \rho^2 (w_x^2+w_y^2)\\
    c_2 &=4 \rho (M w_x -N w_y)\\
    c_3 &=8\rho^2\\
    c_4 &=M^2+N^2-8 \rho^2
    \end{aligned}\nonumber
    \right.
    \label{Eq:a1-a4}
\end{equation}
According to Eq.~(\ref{Eq:RLR6}) and Eq.~(\ref{Eq:RLR9}), we have
\begin{align}
\gamma - \alpha = \theta_f - \dfrac{\pi}{2} + 2 \arccos \dfrac{-\left[X_{T0}+\rho w_x(\dfrac{\pi}{2}- \theta_{ f}+2\beta+2n\pi)\right] - \rho \sin \theta_f + \rho }{4\rho \sin (-\dfrac{\beta}{2})}
\label{Eq:gamma-alpha}
\end{align}
Combining Eq.~(\ref{Eq:gamma-alpha}) with the first equation in Eq.~(\ref{EQ:theta0_thetaf}), we eventually have Eq.~(\ref{EQ:RLR_alpha}) and Eq.~(\ref{EQ:RLR_gamma}).

%%%%%%%%%%%%%%%%%%%%%%%%%%%%%%%%%%%%%%%%%%%%%%%%%%%
Let us now proceed to prove the second statement.
The time required for the Dubins airplane to travel from its initial state $(x_0,y_0,\theta_0)$ to the desired terminal position $(x_f,y_f)$ is equal to the time required for the target to reach the same terminal position in the air-relative frame, indicating
\begin{align}
	\rho (\alpha + \beta + \gamma) =\dfrac{\sqrt{(x_f -X_{T0})^2 + (y_f - Y_{T0})^2}}{\sqrt{(-w_x)^2 + (-w_y)^2}}
	\label{Eq:RLR0_2}
\end{align}
According to Fig.~\ref{Fig:geometry_LRL}, we have 
\begin{equation}
   \left\{
    \begin{aligned}
    &\dfrac{\pi}{2} + \alpha - \beta + \gamma +2n\pi= \theta_f \\
    &\boldsymbol{c}_0^l + 2\rho 
    \left[
    \begin{array}{c}
    \cos\alpha\\
    \sin\alpha
    \end{array}
    \right]=
    \boldsymbol{c}_f^l + 2\rho 
    \left[
    \begin{array}{c}
    \cos(\theta_f - \gamma - \dfrac{\pi}{2})\\
    \sin(\theta_f - \gamma - \dfrac{\pi}{2})
    \end{array}
    \right]
    \end{aligned}
    \label{EQ:theta0_thetaf_2}
    \right.
\end{equation}
where the positive integer $n$ takes values so that  $(\dfrac{\pi}{2} + \alpha - \beta + \gamma +2n\pi) \in [0,2\pi)$, and $\boldsymbol{c}_0^l$ and $\boldsymbol{c}_f^l$ are the centers of the first arc and the third arc, respectively.
Rearranging the second equation in Eq.~(\ref{EQ:theta0_thetaf_2}) yields
\begin{equation}
	\left\{
	\begin{aligned}
		&2\rho [\sin (\alpha + \dfrac{\pi}{2}) - \sin(\theta_f - \gamma)] = x_f - \rho \sin \theta_{f} + \rho \\
        &2\rho [\cos (\alpha +\dfrac{\pi}{2}) - \cos (\theta_f - \gamma)]= -y_f - \rho \cos \theta_f 
	\end{aligned}
	\label{Eq:RLR5_2}
	\right.
\end{equation} 
Combining the first equation in Eq.~(\ref{EQ:theta0_thetaf_2}) and Eq.~(\ref{Eq:RLR5_2}), we have
\begin{equation}
	\left\{
	\begin{aligned}
		& 4\rho \cos \left(\dfrac{\alpha+\pi/2 +\theta_f - \gamma}{2}\right) \sin (\dfrac{\beta}{2})  =  x_f - \rho \sin \theta_{f} + \rho \\
   & -4\rho \sin \left(\dfrac{\alpha+\pi/2  +\theta_f - \gamma}{2}\right) \sin (\dfrac{\beta}{2})  = -y_f - \rho \cos \theta_f 
	\end{aligned}
	\label{Eq:RLR6_2}
	\right.
\end{equation} 
Squaring both sides of each equation in Eq.~(\ref{Eq:RLR6_2}) and adding the results yields
\begin{align}
16 \rho^2 \sin^2 \dfrac{\beta}{2} = \left[ x_f - \rho \sin \theta_{f} + \rho \right]^2 + \left[ y_f + \rho \cos \theta_f \right]^2
\label{Eq:RLR7_2}
\end{align}
Substituting Eq.~(\ref{EQ:theta0_thetaf_2}) into Eq.~(\ref{Eq:RLR0_2}), we have
\begin{align}
\rho(\theta_f-\dfrac{\pi}{2} - 2n\pi) + 2 \rho \beta = \dfrac{\sqrt{(x_f -X_{T0})^2 + (y_f -Y_{T0})^2}}{\sqrt{w_x^2 + w_y^2}}
\label{Eq:RLR8_2}
\end{align}
Notice that
\begin{align}
    -w_y(x_f-X_{T0}) =-w_x(y_f-Y_{T0})
    \label{Eq:RLR4_2}
\end{align}
Combining Eq.~(\ref{Eq:RLR4_2}) and Eq.~(\ref{Eq:RLR8_2}) leads to
\begin{equation}
\left\{
\begin{aligned}
x_f=X_{T0}+\rho w_x(\theta_{ f}-\dfrac{\pi}{2}+2\beta -2n\pi)\\
y_f = Y_{T0} +\rho w_y(\theta_{ f}-\dfrac{\pi}{2}+2\beta -2n\pi)
    \label{Eq:RLR9_2}
\end{aligned}
\right.
\end{equation}
Substituting Eq.~(\ref{Eq:RLR4_2}) and Eq.~(\ref{Eq:RLR9_2}) into Eq.~(\ref{Eq:RLR7_2}) to eliminate $x_f$ and $y_f$, we obtain a transcendental equation, i.e.,
\begin{align}
    16\rho^2\sin^2\dfrac{\beta}{2}=(P+2\rho w_x\beta)^2+(Q-2\rho w_y\beta)^2
\label{Eq:RLR10_2}
\end{align}
where
\begin{equation}
    \left\{
    \begin{aligned}
    P &=X_{T0}+\rho w_x(\theta_{ f}-\dfrac{\pi}{2} -2n\pi)
    -\rho \sin \theta_f + \rho \\
    Q &=-Y_{T0} - \rho w_y (\theta_{ f}-\dfrac{\pi}{2} -2n\pi) - \rho \cos \theta_f 
    \end{aligned}\nonumber
    \right.
\end{equation}
Rearranging Eq.~(\ref{Eq:RLR10_2}), we have 
\begin{align}
    d_1\beta^2 + d_2 \beta +d_3\cos \beta +d_4 =0
\end{align}
where $d_1$--$d_4$ are given as
\begin{equation}
    \left\{
    \begin{aligned}
    d_1 &=4 \rho^2 (w_x^2+w_y^2)\\
    d_2 &=4 \rho (P w_x -Q w_y)\\
    d_3 &=8\rho^2\\
    d_4 &=P^2+Q^2-8 \rho^2
    \end{aligned}\nonumber
    \right.
    \label{Eq:a1-a4_2}
\end{equation}
According to Eq.~(\ref{Eq:RLR6_2}) and Eq.~(\ref{Eq:RLR9_2}), we have
\begin{align}
\alpha - \gamma = -\theta_f - \dfrac{\pi}{2} + 2 \arccos \dfrac{-\left[X_{T0}+\rho w_x(\dfrac{\pi}{2}- \theta_{ f}+2\beta+2n\pi)\right] - \rho \sin \theta_f + \rho }{4\rho \sin (\dfrac{\beta}{2})}
\label{Eq:gamma-alpha_2}
\end{align}
Combining Eq.~(\ref{Eq:gamma-alpha_2}) with the first equation in Eq.~(\ref{EQ:theta0_thetaf_2}), we eventually have Eq.~(\ref{EQ:LRL_alpha}) and Eq.~(\ref{EQ:LRL_gamma}).
This completes the proof.

%%%%%%%%%%%%%%%%%%%%%%%%%%%%%%%%%%%%%%%%%%%%%%%%%%%%
\subsection*{A.3 Expressions of $e_1$--$e_3$ and $f_1$--$f_5$ for Lemma \ref{LE:Zheng}}

According to Appendix B in \cite{ZHENG2021}, if the solution path is of $RSR$, the expressions of $e_1$--$e_3$ in Eq.~(\ref{EQ:RSR}) are given as
\begin{equation}
    \left\{
        \begin{aligned}
        e_1&=(\rho\sin\theta_f-\rho)w_y+\left( \dfrac{w_y}{w_x}X_{T0}-Y_{T0}+\rho\cos\theta_f
        \right)w_x\\
        e_2&=-\rho+\rho\sin\theta_f+\left[ \dfrac{X_{T0}}{w_x}-\rho(\dfrac{\pi}{2}-\theta_f+2n\pi)
        \right]w_x\\
        e_3&=\frac{w_y}{w_x}X_{T0}-Y_{T0}+\rho\cos\theta_f-\left[ \dfrac{X_{T0}}{w_x}-\rho(\dfrac{\pi}{2}-\theta_f+2n\pi)
        \right]w_y
        \end{aligned}
    \right.
    \nonumber
\end{equation}
where $n=0$ if $\theta_f<\pi/2$, and $n=1$ if $\theta_f \geq \pi/2$. If the solution path is of $LSL$, the expressions of $e_1$--$e_3$ in Eq.~(\ref{EQ:RSR}) are given as
\begin{equation}
    \left\{
        \begin{aligned}
        e_1&=-(\rho\sin\theta_f-\rho)w_y+\left( \dfrac{w_y}{w_x}X_{T0}-Y_{T0}-\rho\cos\theta_f
        \right)w_x\\
        e_2&=\rho-\rho\sin\theta_f+\left[ \dfrac{X_{T0}}{w_x}-\rho(\theta_f-\dfrac{\pi}{2}+2n\pi)
        \right]w_x\\
        e_3&=\frac{w_y}{w_x}X_{T0}-Y_{T0}-\rho\cos\theta_f-\left[ \dfrac{X_{T0}}{w_x}-\rho(\theta_f-\dfrac{\pi}{2}+2n\pi)
        \right]w_y
        \end{aligned}
    \right.
    \nonumber
\end{equation}
where $n=0$ if $\theta_f\geq\pi/2$, and $n=1$ if $\theta_f < \pi/2$.

If the solution path is of $RSL$, the expressions of $f_1$--$f_5$ in Eq.~(\ref{EQ:RSL}) are given as
\begin{equation}
    \left\{
        \begin{aligned}
        f_1&=2\rho-(\rho+\rho\sin\theta_f )w_y+\left( -\rho\cos\theta_f+\dfrac{w_y}{w_x}X_{T0}-Y_{T0} \right)w_x\\
        f_2&=2\rho w_y+\left[ \dfrac{X_{T0}}{w_x}-\rho(\dfrac{\pi}{2}+\theta_f+2n\pi) \right]w_x-\rho-\rho \sin\theta_f\\
        f_3&=2\rho w_x-\left[ \dfrac{X_{T0}}{w_x}-\rho(\dfrac{\pi}{2}+\theta_f+2n\pi)\right]w_y-\rho\cos\theta_f+\dfrac{w_y}{w_x}X_{T0}-Y_{T0}\\
        f_4&=-2\rho w_y\\
        f_5&=2\rho w_x
        \end{aligned}
    \right.
    \nonumber
\end{equation}
where
\begin{equation}
    n=\left\{
    \begin{aligned}
        0, \quad \quad &\beta<\pi/2 \ \rm{and} \ \beta<\theta_f\\
        1, \quad \quad &(\beta-\pi/2)(\beta-\theta_f)<0\\
        2, \quad \quad &\beta>\pi/2 \ \rm{and}  \ \beta>\theta_f
    \end{aligned}
    \right.
\end{equation}
If the solution path is of $LSR$, the expressions of $f_1$--$f_5$ in Eq.~(\ref{EQ:RSL}) are given as
\begin{equation}
    \left\{
        \begin{aligned}
        f_1&=-2\rho+(\rho+\rho\sin\theta_f )w_y+\left( \rho\cos\theta_f+\dfrac{w_y}{w_x}X_{T0}-Y_{T0} \right)w_x\\
        f_2&=-2\rho w_y+\left[ \dfrac{X_{T0}}{w_x}+ \rho(\dfrac{\pi}{2}+\theta_f+2n\pi) \right]w_x+\rho+\rho \sin\theta_f\\
        f_3&=-2\rho w_x-\left[ \dfrac{X_{T0}}{w_x}+\rho(\dfrac{\pi}{2}+\theta_f+2n\pi)\right]w_y+ \rho\cos\theta_f+\dfrac{w_y}{w_x}X_{T0}-Y_{T0}\\
        f_4&=-2\rho w_y\\
        f_5&=2\rho w_x
        \end{aligned}
    \right.
    \nonumber
\end{equation}
where
\begin{equation}
    n=\left\{
    \begin{aligned}
        0, \quad \quad &\beta>\pi/2 \ \rm{and} \ \beta>\theta_f\\
        1, \quad \quad &(\beta-\pi/2)(\beta-\theta_f)<0\\
        2, \quad \quad &\beta<\pi/2 \ \rm{and}  \ \beta<\theta_f
    \end{aligned}
    \right.
\end{equation}

%%%%%%%%%%%%%%%%%%%%%%%%%%%%%%%%%%%%%%%%%%%%%%%%%%%%%%%%%
\section{Improved Bisection Method for Solving Transcendental Equations}
\label{Appendix:B}

It is obvious that Eq.~(\ref{EQ:RLR}) and Eq.~(\ref{EQ:LRL}) share the same form. We will propose an improved bisection method for identifying all the roots of the transcendental equation in Eq.~(\ref{EQ:RLR}).

Let us define a function
\begin{align}
	G(\beta) = c_1\beta^2 + c_2 \beta +c_3\cos \beta +c_4 
	\label{EQ:G(alpha)}
\end{align}
where $\beta\in[0,2\pi)$, and $c_1$--$c_4$ are constants.
Differentiate Eq.~(\ref{EQ:G(alpha)}) once and twice, we have
\begin{equation}
    \left\{
    \begin{aligned}
        G'(\beta) &=2c_1 \beta +c_2 -c_3 \sin \beta \\
        G''(\beta) &=2c_1 -c_3 \cos \beta
    \end{aligned}
    \right.
    \label{EQ:G'}
\end{equation}
From the second equation in Eq.~(\ref{EQ:G'}), it can be observed that $G''(\beta)$ has at most two real roots over the interval $[0,2\pi)$. Let us denote these real roots by $\beta_1,...,\beta_n$, where $n\leq 2$ represents the number of roots. Without loss of generality, we assume $\beta_0<\beta_1<...<\beta_n<\beta_{n+1}$, where $\beta_0=0$, and $\beta_{n+1}=2\pi$.

For each interval $[\beta_{i},\beta_{i+1})$, where $i=0,1,...,n$, it is obvious that the first order differentiation $G'(\beta)$ is monotonic. As a result, the number of roots for $G'(\beta)$ within the interval $[\beta_{i},\beta_{i+1})$ is determined by the signs of the boundary values $G'(\beta_{i})$ and $G'(\beta_{i+1})$, indicating
\begin{equation}
    N[G'(\beta)=0 \mid \beta \in [\beta_{i},\beta_{i+1})]=
    \left\{
    \begin{aligned}
        &0, \quad &&\rm{if} \ G'(\beta_i)\times G'(\beta_{i+1}) > 0 \\
        &1, \quad &&\rm{if} \ G'(\beta_i)\times G'(\beta_{i+1}) < 0 \ \rm{or} \  G'(\beta_i)=0
    \end{aligned}
    \right.
    \label{EQ:N(G)}
\end{equation}
If $N[G'(\beta)=0 \mid \beta \in [\beta_{i},\beta_{i+1})] \neq 0$, we use a simple bisection method to determine the unique root. Let us denote by 
\begin{align}
	\gamma = \pmb{B}[G'(\beta),\beta_i,\beta_{i+1}]
	\nonumber
\end{align}
the unique root of the function $G'(\beta)$ within the interval $[\beta_i,\beta_{i+1})$ obtained via the bisection method. Then, denote by $Z_1$ the set of roots for $G'(\beta)$, and it can be readily obtained, i.e.,
$$Z_1=\{\gamma_1,...,\gamma_{n'}\}$$
where $n'$ denote the total number of roots for $G'(\beta)$ over the interval $[0,2\pi)$.

Denote by $Z_2$ the set of all possible roots for $G(\beta)$, and it can also be obtained similarly. The entire process is summarized in Algorithm \ref{algorithm_1}.

\begin{algorithm} 
	\caption{improved bisection method}   % 标题
	 \label{algorithm_1}       % 用来引用
	\begin{algorithmic}[1] % 加上 [1] 表示有序号
    \State $n'\gets 0$
    \State $Z_1\gets \emptyset$
    \State $Z_2\gets \emptyset$
    \For {$i=0$ to $n$}
        \If{$G'(\beta_i) = 0$} 
            \State $Z_1 \gets Z_1 \cup \{\beta_i\}$
            \State $n'\gets n'+1$
        \ElsIf {$G'(\beta_i)\times G'(\beta_{i+1}) < 0$}
            \State $ \gamma \gets B[G'(\beta),\beta_i,\beta_{i+1}]$
		\State $Z_1 \gets Z_1 \cup \{\gamma\}$
            \State $n'\gets n'+1$
        \EndIf
    \EndFor
    \State Sort the elements in set $Z_1$ so that $\gamma_0<\gamma_1<...<\gamma_{n'}<\gamma_{n'+1}$.
    \For {$i=0$ to $n'$}
        \If{$G(\gamma_i) = 0$} 
            \State $Z_2 \gets Z_2 \cup \{\gamma_i\}$.
        \ElsIf {$G(\gamma_i)\times G(\gamma_{i+1}) < 0$}
            \State $ \gamma \gets B[G(\alpha),\gamma_i,\gamma_{i+1}]$.
		\State $Z_2 \gets Z_2 \cup \{\gamma\}$.
        \EndIf
    \EndFor
   \end{algorithmic} 
\end{algorithm}

Upon reviewing Algorithm \ref{algorithm_1},  it is evident that the iteration counts $n$ and $n'$ are finite. Additionally, each iteration of the loop requires at most one finite-time bisection.
This ensures that the set of solutions can be obtained efficiently and robustly.
Testing on a laptop with an Intel Core i5-10400 CPU @ 2.90 GHz shows that, given any function $G(\beta)$, the set $Z_2$ can be obtained in around $10^{-5}$ seconds, thereby comfortably meeting the real-time requirement.

\end{document}